\documentclass[10pt]{elsarticle}

\nonstopmode

\usepackage{subcaption}

\usepackage[utf8]{inputenc}
\usepackage{graphicx}
\usepackage{color}
\usepackage{verbatim}
\usepackage{hyperref}
\usepackage{url}

\usepackage{hypernat}
\usepackage{mathpazo}
\usepackage{cleveref}
\usepackage{amsmath}
\usepackage{amsfonts}
\usepackage{amssymb}
\newtheorem{defin}{Definition}
\newtheorem{thm}{Theorem}
\newtheorem{cond}{Condition}
\newtheorem{lem}{Lemma}

\newtheorem{rem}{Remark}

\numberwithin{equation}{section}

\usepackage{setspace}

\usepackage{stmaryrd}

\usepackage{tikz}
\usetikzlibrary{positioning}
\usetikzlibrary{calc}
\usetikzlibrary{chains}
\usetikzlibrary{shapes}
\usetikzlibrary{arrows}

\pgfdeclarelayer{background}
\pgfdeclarelayer{foreground}
\pgfsetlayers{background,main,foreground}

\usepackage{floatrow}
\newfloatcommand{capbtabbox}{table}[][\FBwidth]

\usepackage[noend]{algorithmic}
\usepackage{algorithm}
\algsetup{indent=2em}

\usepackage{multirow,hhline}

\usepackage[margin=10pt,font=small,labelfont=bf,labelsep=period]{caption}

\newcommand{\doi}[1]{\href{http://dx.doi.org/#1}{doi: #1}}

\usepackage{geometry}
\newcommand{\bx}{\boldsymbol {x}}
\newcommand{\bt}{\boldsymbol {t}}

\newcommand{\bs}{\boldsymbol {s}}

\newcommand{\bc}{\boldsymbol {c}}
\newcommand{\bmm}{\boldsymbol {m}}
\newcommand{\bn}{\boldsymbol {n}}
\newcommand{\bbR}{\mathbb{R}}
\newcommand{\cC}{\mathcal{C}}
\newcommand{\cS}{\mathcal{S}}
\newcommand{\cD}{\mathcal{D}}
\newcommand{\cO}{\mathcal{O}}

\newcommand{\pfmm}{p_\text{FMM}}
\newcommand{\pqbx}{p_\text{QBX}}
\newcommand{\pqbxl}{p_\text{QBX}^\ell}
\newcommand{\padd}{p_\text{add}}
\newcommand{\pfmml}{p_\text{FMM}^{\ell}}

\newcommand{\tgtquant}[1]{{#1}^\text{dens}}
\newcommand{\quadquant}[1]{#1}
\newcommand{\qtgt}{\tgtquant{q}}

\let\fourquad=\qquad
\renewcommand{\qquad}{\quadquant{q}} 
\newcommand{\sigquad}{\quadquant{\sigma}}
\newcommand{\squad}{\quadquant{\boldsymbol{s}}%
}

\newcommand{\stgt}{\tgtquant{\boldsymbol{s}}%
}
\newcommand{\ntgt}{\boldsymbol{n}}

\newcommand\qline{\begin{center}\line(1,0){300}\end{center}}

\newcommand{\cB}{\mathcal{B}}

\let\epsilon=\varepsilon

\def\marklegendre#1{
    \foreach \i / \ln in {
       0/0.02544604,  1/0.12923441, 2/0.29707742, 3/0.5,
       4/0.70292258, 5/0.87076559, 6/0.97455396}
    { coordinate [pos=\ln] (#1\i) }
  }

\pgfdeclarelayer{background}
\pgfdeclarelayer{grid}
\pgfdeclarelayer{foreground}
\pgfsetlayers{background,grid,main,foreground}

\newenvironment{tabbingalg}{
  \begin{minipage}{\textwidth}
  \small
  \begin{tabbing}
  ......\=......\=......\=......\=\kill
}{
  \end{tabbing}
  \end{minipage}\\
}

\begin{document}


\title{%
  Fast algorithms for Quadrature by Expansion I:
  Globally~valid~expansions
}

\author[yale]{Manas Rachh \corref{MR}}%
\ead{manas.rachh@yale.edu}
\author[uiuc]{Andreas Kl\"ockner}%
\ead{andreask@illinois.edu}
\author[nyu]{Michael O'Neil}%
  \ead{oneil@cims.nyu.edu}
\cortext[MR]{Corresponding author}
\address[yale]{Applied Mathematics Program, Yale University, 51
  Prospect St, New Haven, CT 06511}
\address[uiuc]{Department of Computer Science,
  University of Illinois at Urbana-Champaign, 201 North Goodwin Ave,
  Urbana, IL 61801}
\address[nyu]{Courant Institute and Tandon School of Engineering,
New York University, New York, NY}

\begin{abstract}
  The use of integral equation methods for the efficient numerical
  solution of PDE boundary value problems requires two main tools:
  quadrature rules for the evaluation of layer potential integral
  operators with singular kernels, and fast algorithms for
  solving the resulting dense linear systems.
  Classically, these tools were developed separately.
  In this work, we present a unified numerical scheme based on
  coupling {\em Quadrature by Expansion}, a recent quadrature method, to a
  customized Fast Multipole Method (FMM) for the Helmholtz equation
  in two dimensions.
  The method allows the evaluation of layer potentials in linear-time
  complexity, anywhere in space, with a uniform, user-chosen level of
  accuracy as a black-box computational method.

  Providing this capability requires geometric
  and algorithmic considerations beyond the needs of standard FMMs
  as well as careful
  consideration of the accuracy of multipole translations.
  We illustrate the speed and accuracy of our method with
  various numerical examples.
\end{abstract}


\maketitle


\onehalfspacing

\section{Introduction}\label{sec_intro}

Reformulating the partial differential equations (PDEs) of classical
mathematical physics in integral form and then discretizing the
resulting integral equation affords several analytic
and computational
advantages over direct discretizations of the differential operator.
For example, integral representations of the solution to exterior
boundary value problems inherently capture the correct decay
properties at infinity via the use of the Green's function for the PDE\@.
Furthermore, integral equation formulations are often able to reduce
volume discretizations to boundary discretizations, yielding an immediate
reduction in computational and storage complexity. Lastly, integral equation formulations
frequently reflect the natural conditioning of the underlying physical
problem, i.e., well-conditioned physical problems yield
well-conditioned integral equations.
For concreteness, we consider the method in the setting of the exterior
Dirichlet problem for the
Helmholtz equation in two dimensions in this contribution.
We note however that the method generalizes rather straightforwardly
to higher dimensions and different kernels and layer potentials.
Consider the boundary value problem~\cite{guenther_lee}:
  \begin{align}
    \left( \triangle + \omega^2 \right) u &= 0 \quad \text{in } \mathbb
    R^2 \setminus \Omega, \label{eq:helm-eq} \\
    u &= f \quad \text{on } \partial \Omega ,
    \label{eq_diri} \\
  \lim_{r\to\infty}
  r^{1/2}
  \left(\frac{\partial}{\partial r} - i\omega \right) u &= 0,
  \label{eq:sommerfeld-rad}
  \end{align}
where $\Omega \subset \mathbb R^2$ is a closed, bounded region with
smooth boundary $\Gamma = \partial \Omega$.
Equation \eqref{eq:helm-eq} is the
Helmholtz equation in $\bbR^{2} \setminus \Omega$,
equation \eqref{eq_diri} enforces Dirichlet boundary conditions,
and equation \eqref{eq:sommerfeld-rad} is the {\em Sommerfeld radiation
condition} which ensures that the solution $u$ is a {\em radiating solution}.
This boundary-value problem
can be reformulated in integral form by representing the solution
$u$ as a {\em combined-field potential} \cite{brakhage_uber_1965,panich_solubility_1965}:
\begin{equation}\label{eq_dlp}
\begin{aligned}
u(\bx) &= \mathcal D[\sigma] (\bx) + i\omega \mathcal S[\sigma](\bx)\\
&= \int_\Gamma \frac{\partial G}{\partial n_{x'}}(\bx,\bx')
  \, \sigma(\bx') \, ds(\bx') + i\omega
\int_\Gamma G(\bx,\bx')
  \, \sigma(\bx') \, ds(\bx')
\end{aligned}
\end{equation}
with $G$ the Green's function
\begin{equation}
  G(\bx,\bx') = \frac{i}{4} H^{(1)}_0(\omega|\bx-\bx'|),
  \label{eq:helmholtz-ker}
\end{equation}
where $H^{(1)}_0$ denotes the zeroth-order Hankel function of the
first kind \cite{abramowitz_1965}, and $\partial/\partial n_{x'} = \bn_{x'} \cdot
\nabla_{x'}$ with $\bn_{x'}$ the unit normal vector at
$\bx'$ pointing out of $\Omega$.
Here, unless otherwise specified, $|\cdot|$ denotes the $\ell_2$-norm
when applied to a vector,
$\omega \in\mathbb C$
with $\operatorname{Im} \omega \ge 0$, and $\sigma$ is an unknown
density defined on the boundary $\Gamma$.
This Green's function also satisfies the Sommerfeld radiation  condition
\[
  \lim_{r\to\infty}
  r^{1/2}
  \left(\frac{\partial}{\partial r} - i\omega \right)
H^{(1)}_0(\omega  r)  = 0.
\]
Using this representation of $u$ along with the so-called
jump relations (cf. \eqref{eq_oneside})
and enforcing the Dirichlet boundary condition in~\eqref{eq_diri},
we obtain the following integral
equation along $\Gamma$:
\begin{equation}
\label{eq_extdir}
\frac{1}{2} \sigma  + \mathcal D^* [\sigma] + i\omega \mathcal S^* [\sigma]
= f.
\end{equation}
The operators $\mathcal D^*$ and
$\cS^*$, as  maps from
$\Gamma \to \Gamma$, are merely $\cD$ and $\cS$
interpreted as principal value and improper integrals, respectively. The
operator $\cS$ is known as the {\em single-layer} potential, and $\cD$
as the {\em double-layer} potential.
The efficient solution
of~\eqref{eq_extdir} requires the development of
several numerical tools, including quadrature methods
for layer potentials with kernels of varying degrees of singularity
and, in the setting of iterative methods such as GMRES,
asymptotically fast algorithms for computing matrix-vector
products with the (notionally dense) matrices resulting from the
discretization of \eqref{eq_extdir}.

Historically, these two numerical tools -- quadrature for singular
functions and fast algorithms for applying discretized integral
operators -- have been treated separately.
Quadrature methods for singular functions include product
integration~\cite{helsing_2011}, generalized Gaussian quadrature
rules~\cite{bremer,yarvin_1998},
singularity subtraction~\cite{jarvenpaa_2003}, and many
others.
See~\cite{hao_2014} for a recent overview of existing methods.
The problem of singular quadrature for layer potentials is, in
essence, a local one, in the sense that
when the {\em target} and {\em source}, $\bx$ and $\bx'$,
respectively, are well-separated, conventional quadrature rules (such as composite Gaussian)
may be applied effectively. As such, the treatment of the singularity
may be constrained to the near-field of each piece of the
geometry $\Gamma$.
Such quadrature rules are generally straightforward to couple
with fast algorithms
such as FMMs because of the natural
separation of near-field and far-field calculations.

A somewhat nuanced sub-problem in terms of quadrature
appears when the integrand is not singular, but only nearly so.
This occurs when targets are
located \emph{near} but not \emph{on} the boundary~$\Gamma$.
In this sense, targets close to the source geometry present a
different challenge than on-surface targets: the corresponding integrals
are computable, e.g. by adaptive quadrature, however maintaining
efficiency has proven nontrivial for previous methods.
Quadrature by Expansion (QBX)~\cite{qbx,epstein} originated as an extension of
a scheme for nearby evaluation~\cite{barnett_evaluation_2012},
and thus provides a means to handle this difficulty.

In this work, we present a unified algorithm which efficiently embeds
Quadrature by Expansion inside a fast
multipole method for the two-dimensional Helmholtz equation.
The algorithm evaluates both on-surface layer potentials as well
as potentials at points arbitrarily near the boundary $\Gamma$
in a single computation.

QBX was initially described \cite{qbx} in terms of the underlying
analytical idea along with initial theoretical and numerical insight.
No fast algorithm was applied, the relevant examples required
 $\mathcal O(n^2)$ operations
for potential evaluation along a boundary discretized using
$n$ points.
Subsequent work on QBX
in~\cite{epstein,barnett_evaluation_2012,af2016estimation}
is largely foundational, and
investigates the rate of asymptotic and convergent approximations of
such potential expansions.
The goal of this paper, on the other hand, is to make QBX a viable
numerical algorithm for the solution to boundary integral
equations in two dimensions.
The derivation of an algorithm for three dimensions is
straightforward using the analogous hierarchical data structures as are described
in this work.
In order for this
to be accomplished, QBX has to be coupled with a fast algorithm
such as an FMM.
There are several additional considerations that need
to be addressed: the effect of adaptive geometry discretizations,
orders of multipole and local expansions, and the effective \emph{radii of
accuracy} for translated local expansions.
An early version of a fast QBX algorithm was used
in~\cite{oneil_imped}, but details regarding these topics were not
provided.

The paper is organized as follows: In Section~\ref{sec:background} we
briefly review the details of Quadrature by Expansion for Helmholtz
layer potentials and give an informal description of an
accelerated algorithm using a standard FMM.
In Section~\ref{sec:accuracy}, we discuss elements of the underlying
geometry discretization that will affect the resulting accuracy of the
scheme.
Sections~\ref{sec:area-query} and~\ref{sec:source-refine} discuss
additional data structure elements, and how these
features adaptively refine the geometry and layer potential densities
to ensure numerical accuracy.
In Section~\ref{sec:volume-eval}, we discuss the details of using the
global QBX scheme to evaluate potentials off (but near) the surface.
A detailed description of the QBX scheme embedded into an FMM is given in
Sections~\ref{sec:fmm}, followed by numerical examples in
Section~\ref{sec:results}.
Lastly, in Section~\ref{sec:conclusions}, the conclusion,
we discuss drawbacks of the
algorithm of this paper, as well as describe extensions to three
dimensions and other PDEs.

\section{Background material}
\label{sec:background}
In this section we give an overview of the existing QBX scheme,
assumptions on the underlying geometry discretizations and the layer
potential densities (i.e $\sigma$), and then provide an informal
description of the algorithm.  It is important to make a distinction
between the on-surface value of the layer potentials
$\mathcal S[\sigma]$ and~$\mathcal D[\sigma]$ (as an improper or
principal-value integral) and their one-sided limits.  The one-sided
limits of these layer potentials are \cite{kress_1999}
\begin{equation}\label{eq_oneside}
  \begin{aligned}
    \lim_{s \to 0^\pm} \mathcal S[\sigma](\bx + s  \bn_x) &=
    \mathcal S^*[\sigma](\bx), \\
    \lim_{s \to 0^\pm} \mathcal D[\sigma](\bx + s \bn_x) &=
    \mp\frac{1}{2}\sigma(\bx) + \mathcal D^*[\sigma](\bx),
  \end{aligned}
\end{equation}
where we have assumed that the point $\bx$ is located on the curve
$\Gamma$.
Note that the potential $\cS[\sigma]$ is continuous across the
boundary $\Gamma$, and in both cases, the potentials are smooth up to
the boundary $\Gamma$ with well-defined limits.
By~$\cS^*$ and~$\cD^*$ we denote the on-surface restriction of the
operators~$\cS$ and~$\cD$~\cite{kress_1999, colton_kress}.

Quadrature by Expansion computes the one-sided limits on the
right-hand side in~\eqref{eq_oneside}.
This is in contrast to classical quadrature schemes for layer
potentials which compute the on-surface principal-part of the
operator, while one-sided limits are obtained using the jump
relations.
By computing the average of two different one-sided QBX calculations,
 as described in~\cite{qbx}, one can use QBX to compute
the on-surface value of the operator.
This approach has certain advantages when coupled with
iterative solvers, but this usage is outside of the scope of the
current paper.

\subsection{Using expansions for quadrature \label{sec:expquad}}

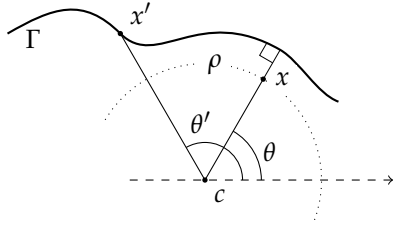
\begin{figure}[!t]
  \begin{center}
  \begin{tikzpicture}
    \coordinate (c) at (0,0) ;
    \path (c) ++(45+75:2.25) coordinate (s);
    \path (c) ++(45+15:1.55) coordinate (t);
    \path (c) ++(45+15:2) coordinate (t-to-curve);

    \path (s) ++(-1.5,0) coordinate (curve-before);
    \path (t) ++(1,-0.3) coordinate (curve-after);
    \draw [thick]
      (curve-before)
      ..controls +(30:0.7) and +(180-45:0.7) ..
      (s)
      node [pos=0.2,anchor=north] {$\Gamma$}
      ..controls +(-45:0.7) and +(45+90+15:1.25) ..
      (t-to-curve)
      ..controls +(45+270+15:0.3) and +(160:0.3) ..
      (curve-after) ;

    \fill (c) circle (1pt);
    \fill (t) circle (1pt);
    \fill (s) circle (1pt);

    \draw [dotted] ++(-20:1.55) arc (-20:45+105:1.55);

    \node at (85:1.55) [fill=white] {$\rho$};
    \draw (c) -- (t) ;
    \draw (c) -- (s) ;
    \node at (c) [anchor=north west] {$c$};
    \node at (s) [anchor=south west] {$x'$};
    \node at (t) [anchor=west,fill=white,xshift=0.5mm,inner sep=1mm] {$x$};
    \draw [dashed,->] (c) ++(-1, 0) -- ++(3.5,0);

    \draw (c) ++(0.5,0) arc (0:45+75:0.5);
    \path (c) ++(45*0.8+75*0.8:0.5)
      node [anchor=south] {$\theta'$};

    \draw (c) ++(0.75,0) arc (0:45+15:0.75) ;
    \path (c) ++(45/2+15/2:0.75)
      node [anchor=west] {$\theta$};

    \draw [very thin] (t) -- (t-to-curve);
    \draw [very thin]
      let
        \p1 = ($ 0.1*(t-to-curve) - 0.1*(c) $),
        \p2 = (-\y1,\x1)
      in
      ($(t-to-curve)!0.1!(c)$) -- ++(\p2) -- ++(\p1) ;

  \end{tikzpicture}
  \caption{Geometric configuration for Graf's addition formula with
    respect to a source $\bx'$ located on $\Gamma$ and a target
    $\bx$ located on or off of $\Gamma$.}
    \label{fig_expans}
  \end{center}
\end{figure}

We next give a precise description of QBX
in the case of evaluating $\mathcal S[\sigma]$ in order to
establish notation.
For a target point
$\bx$ near the boundary $\Gamma$, the computational task is the
numerical evaluation of the integral
\begin{equation}
\begin{aligned}
u(\bx) &= \mathcal S[\sigma](\bx) \\
&= \int_\Gamma G(\bx,\bx')
  \, \sigma(\bx') \, ds(\bx') \label{eq:slp}.
\end{aligned}
\end{equation}
Using the definition of the Green's function $G$
in~\eqref{eq:helmholtz-ker}, and the Graf
addition formula~\cite{olver_nist_2010},
for $\bx$ not on the boundary $\Gamma$,
we can rewrite the potential $u$ as:
\begin{equation}\label{eq_jexp}
u(\bx) = \sum_{\ell = -\infty}^\infty \alpha_\ell \, J_\ell(\omega\rho) \,
e^{-i\ell \theta}.
\end{equation}
The coefficients $\alpha_\ell$ are known as {\em local expansion
coefficients},
Fourier-Bessel coefficients, or coefficients of a $J$-expansion.
They are given explicitly by:
\begin{equation}
\label{eq:graf-coefficient}
\alpha_\ell (\bc) = \frac{i}{4} \int_{\Gamma}
H_\ell^{(1)}( \omega \vert\bx'-\bc\vert) \, e^{i\ell\theta'}
\sigma(\bx') \, ds(\bx'),
\end{equation}
where the polar coordinates of $\bx -\bc = (\rho,\theta)$ and
$\bx'-\bc = (\rho',\theta')$
are with
respect to the expansion center~$\bc$, with $\vert \bx - \bc\vert <
\vert \bx'-\bc\vert$, located
off of $\Gamma$ and usually along the normal to the curve near $\bx$.
Here $H_{\ell}^{(1)}$ denotes the $\ell^\text{th}$ order Hankel function
of the first kind and $J_{\ell}$ is the $\ell^\text{th}$ order
Bessel function of the first kind.
This restriction on the location of $\bc$, relative to $\bx$ and $\bx'$,
is necessary to ensure the validity of Graf's identity. See
Figure~\ref{fig_expans} for a graphical depiction.
The $\omega$-scaled Fourier-Bessel basis satisfies the two-dimensional
Helmholtz equation and thus
is a natural basis in which to expand potentials.

Even though this expansion was constructed about a point $\bc$
away from the boundary $\Gamma$, it can be evaluated at points
$\bx$ near, or even {\em on} the boundary.
In practice, the coefficients $\alpha_\ell$ in such expansions decay
at a rate that depends on the smoothness of $u$ and the distance of
$\bc$ to the boundary. Truncating this series at index $p$ yields
an approximation to the potential, with truncation estimates given in \cite{epstein}:
\begin{equation}\label{eq_phiapprox}
u(\bx) \approx \sum_{\ell = -p}^p \alpha_\ell \, J_\ell(\omega\rho) \,
e^{-i\ell \theta}.
\end{equation}
The approximation
of integrals for the coefficients $\alpha_\ell$ via quadrature,
and subsequent
evaluation of the truncated series~\eqref{eq_phiapprox}, is the
essence of \emph{Quadrature by Expansion}. This procedure
effectively computes the singular or nearly-singular
integral $\cS[\sigma]$ by expanding the
resulting potential about some point $\bc$ off the surface, and
simultaneously, regularizing the integrals corresponding to $\alpha_\ell$.
Similar expansions to~\eqref{eq_phiapprox} exist for various other
kernels and layer potentials.

As calculated above, the coefficients $\alpha_\ell$ will be referred
to as {\em global QBX local expansion coefficients} since they contain
information from all of $\Gamma$.
We will call an application of QBX with such coefficients {\em Global QBX}.
This is in contrast to {\em Local QBX}, to be described in a future
contribution.
The algorithm of this paper is concerned with constructing a fast
and accurate scheme for the calculation of these global coefficients.

In order to derive the a simple scheme for embedding QBX inside an
FMM, it is necessary to place restrictions on the boundary $\Gamma$ and
its discretization. The follow sections outlines our notation for
boundary discretization, as well as the restrictions that we place on it.

\subsection{Geometry discretization}
\label{sec:geom-discr}

We assume that the source curve $\Gamma=\cup_{j=1}^{N} \Gamma_{j}$ from which
we wish to compute the layer potential is discretized by
$N$ piecewise Gauss-Legendre panels $\Gamma_{j}$ with $\qtgt$
Gauss-Legendre points
on each panel on which the density is represented.
We typically discretize panels with $\qtgt=2$, $4$, $8$ and $16$
points. We sometimes term this the `density' discretization.

We note that in addition to this `density' discretization,
we will find it necessary to introduce a `source'
discretization consisting of the same panels, each equipped with
a larger number of quadrature points, so that the integrals involving
density and kernel moments involved in QBX can be
computed accurately. Section~\ref{sec:quaderror} contains the specifics.
We will denote the number of source (or quadrature) points
on a panel $\qquad$. When needed, point values of the density
on the `source' discretization are obtained from the target
discretization by polynomial interpolation.

We will use the following notation to refer to nodes and
other geometric entities:
\begin{itemize}
\item $h_{k} = \int_{\Gamma_k} ds$, the arclength of the $k^\text{th}$
  panel,
\item $\stgt_{j,k}$, the location of the $j^\text{th}$ (`density') node on the
  $k^\text{th}$ panel, discretized with respect to arclength,
\item  $\ntgt_{j,k}$, the outward unit normal to the boundary at the
  $j^\text{th}$  node on the $k^\text{th}$ panel,
\item $\bc_{j,k} = \stgt_{j,k} + \frac{h_{k}}{2} \ntgt_{j,k}$ the expansion
  centers corresponding to node $\stgt_{j,k}$, and
\item $\squad_{j,k}$, the location of the $j^\text{th}$ (`source') node on the
  $k^\text{th}$ panel, discretized with respect to arclength,
\item $\bt_{j}$, for $j=1,2,\ldots n_{t}$ denote the location of the
  targets in the volume at which we wish
  to evaluate the layer potential.

  We note that the targets may reside
  anywhere in $\overline{\mathbb R^2\setminus \Omega}$. They are explicitly
  not restricted to the boundary.
  Nonetheless, some subtlety is required for targets residing
  on or near the source curve, see Section~\ref{sec:volume-eval} for
  details.
\end{itemize}

Along the boundary $\Gamma$, the continuous density $\sigma$ is
discretized (sampled) at each of the Gauss-Legendre nodes on each
panel in accordance with the Nyström method \cite{nystroem_praktische_1930}.
For targets $\bx$ sufficiently far away from the boundary
$\Gamma$, the single-layer potential can then be accurately
computed by a $\qquad$-point Gaussian quadrature rule for smooth functions:
\begin{equation}\label{eq_sigdisc}
\begin{aligned}
u(\bx) &= \mathcal S[\sigma](\bx) \\
&\approx \frac{i}{4} \sum_{k=1}^N \sum_{j=1}^\qquad w_{j,k} \, H_0^{(1)}(\omega\vert
\bx - \squad_{j,k}\vert) \, \sigma(\squad_{j,k}).
\end{aligned}
\end{equation}
If the standard $\qquad$-point
Gauss-Legendre quadrature weights on the interval $[-1,1]$ are given
by $\quadquant{w}_j$, then $\quadquant{w}_{j,k} = h_k \quadquant{w}_j/2$
since it is assumed the nodes
$\squad_{j,k}$ are sampled with respect to arclength.

We now turn to the topic of estimating the error inherent in
QBX as a numerical methdo, assuming a piecewise Gauss-Legendre discretization of
$\sigma$ along $\Gamma$.


\subsection{Estimating the error in QBX}
\label{sec_error}

The error in standard quadrature rules (for smooth functions)
is usually determined by a single
parameter, the \emph{order} of the quadrature.
The error estimates for QBX are slightly more complicated, and consist
of two components: the truncation error (related to the decay of the local expansion
used) and the quadrature error (i.e.\ error in computing the coefficients of
the local expansion).

As shown in \cite{qbx}, the combined truncation and quadrature error
in QBX can be estimated as follows, based on \cite[eqn. (2.7.12)]{davis_1984}:

\begin{thm} \label{thm:qbx-orig}
Suppose that $\Gamma$ is a smooth, bounded curve embedded in $\bbR^2$,
and that $B_r(\bc) \cap \Gamma = \emptyset$, where $B_r(\bc)$ is
the open ball of radius $r>0$ centered at $\bc$.  Let $\Gamma$ be
divided into $M$ panels, each of length $h$, and
let $Q_q(f)$ denote the Gauss-Legendre quadrature approximation to the
integral $f$
using $\qquad$ points.
For $0<\beta<1$, there are constants
$C_{p,\beta,\Gamma}$ and $\tilde C_{p,\qquad,\beta,\Gamma}$ so that if $\sigma$ lies in
the Hölder space $\cC^{p,\beta}(\Gamma)\cap \cC^{2\qquad,\beta}(\Gamma)$,
then
\begin{equation}
  \left\vert \mathcal S[\sigma](\bx) -
  \sum_{l=-p}^{p} Q_q(\alpha_\ell) \,
  J_\ell(\omega|\bx-\bc|) \, e^{-i\ell\theta}\right\vert
  \leq
    \underbrace{
      C_{p,\beta,\Gamma} \; r^{p+1}
      \|\sigma\|_{\cC^{p,\beta}(\Gamma)}
    }_{\text{Truncation error}}
    +
    \underbrace{
      \tilde C_{p,\qquad,\beta,\Gamma} \left( \frac{h}{4r} \right)^{2\qquad}
      \|\sigma\|_{\cC^{2\qquad,\beta}(\Gamma)}
    }_{\text{Quadrature error}},
  \label{eq:underlying-quad-estimate-v1}
\end{equation}
where $\mathcal S[\sigma]$ is the single-layer potential defined
in~\eqref{eq_dlp}, and $\theta$ is as in Figure~\ref{fig_expans}.
\end{thm}
\begin{rem}
The quadrature error contains a factor of $(h/4r)^{2\qquad}$. Thus,
as long as $r>h/4$, this factor can be used to
control quadrature error by increasing $\qquad$.
However, increasing $\qquad$ \emph{still} decreases the quadrature error
substantially even if $r\leq h/4$, which the above estimate \emph{fails}
to predict. (We
encounter this situation in the determination of $\qquad$ in
Section~\ref{sec:quaderror}.)
More precise error estimates for the QBX quadrature error
based on characteristic discretization lengths
(i.e. not adaptive discretizations) are discussed
in~\cite{af2016estimation}.
\end{rem}

\subsection{Informal description of the algorithm}

With the method of QBX and the previous error estimates in mind, an
accelerated scheme for the computation of $\cS[\sigma]$ can be derived.
This global, accelerated FMM-based QBX algorithm involves three basic
steps:
\begin{quote}
\begin{enumerate}
\item[\bf Step 1] First, we refine the discretization of the boundary
  $\Gamma$ and the density $\sigma$
  to ensure the validity of the error
  estimate~\eqref{eq:underlying-quad-estimate-v1}.
  For this discretization, we then place one QBX expansion
  center per `density' node.
\item[\bf Step 2] Using a modification of the FMM for the two-dimensional
  Helmholtz equation, we evaluate the potential at all targets,
  and construct a (QBX) local expansion at the
  expansion centers that were placed in Step 1 using translations of
  the FMM-computed Bessel function expansions.
\item[\bf Step 3] Finally, we identify targets that are \emph{close}
  to the boundary (i.e.\ such targets where the underlying Gaussian
  quadrature rule fails to accurately approximate the singular or
  near-singular integral), and re-evaluate the potential using the
  (QBX) local expansions of the expansion center closest to the
  target.
\end{enumerate}
\end{quote}
The remainder of the paper is dedicated to filling in the details
surrounding each of these steps.



\section{Accuracy control for global QBX}
\label{sec:accuracy}

As noted in Theorem~\ref{thm:qbx-orig},
QBX incurs two (additive) error components,
truncation error and quadrature error.
The truncation error stems from using a truncated $J$-expansion
to approximate the (smooth) layer potential, whereas
the quadrature error arises from computing
the expansion coefficients $\alpha_\ell$
using numerical quadrature.
We now discuss some practical sufficient conditions under which
the method yields potentials that are pointwise convergent of order
$p+1$ in the maximum panel length $\max_{k} h_k$
up to a user-supplied precision $\epsilon$.
To achieve this,
we ensure that (a) the assumptions of Theorem~\ref{thm:qbx-orig}
apply, and (b) the quadrature error does not exceed $\epsilon$. Note that
we make no attempt to characterize truncation error beyond high-order
convergence.

\subsection{Controlling truncation error
\label{sec:trunc-error-control}}

Truncation error in QBX is tied to the decay of the coefficients
of the local (Fourier-Bessel) expansion used to evaluate the
potential near the surface. The decay of these coefficients
reflects the smoothness of the expanded potential, and the smoothness
of the potential in turn is controlled by the proximity of any source
geometry. Theorem~\ref{thm:qbx-orig} assures high-order convergence
as long as there is no source geometry on the interior of the
expansion disk. As such, controlling truncation error in QBX is mainly a
geometric matter.

Given our choice of expansion radius
$r=h_k/2$ for all centers associated with sources on panel $k$
(cf. Sec.~\ref{sec:geom-discr}), satisfying the assumptions
of Theorem~\ref{thm:qbx-orig} amounts to satisfying the following
condition.
Here and in the following, $d(\cdot,\cdot)$ denotes the Euclidean
distance function as applied to points and sets of points.

\begin{cond}[Expansion disk undisturbed by sources]
  \label{cond:no-sources-in-disk} To ensure that no source
  geometry interferes with the decay of the coefficients
  of the local expansion, we demand that the distance
  between the source geometry and each center be at least
  the radius of that center's expansion disk, i.e.
  \begin{equation}\label{eq_dcgam}
  d (\bc_{j,k},\Gamma) \geq \frac{h_{k}}{2}
  \end{equation}
  for all nodes $j$ and panels $k$.
  Breaking down the source geometry further, we may equivalently
  enforce
  \begin{equation}\label{eq_dcgam_perpanel}
  d (\bc_{j,k},\Gamma_\ell) \geq \frac{h_{k}}{2}
  \end{equation}
  for all nodes $j$ and panels $k$ and $\ell$.
\end{cond}

Under Condition~\ref{cond:no-sources-in-disk}, the truncated
Fourier-Bessel expansion
\[
\sum_{\ell=-p}^{p} \alpha_{\ell} \,
 J_\ell(\omega|\bt-\bc_{j,k}|) \, e^{-i\ell\theta}
\]
is a $(p+1)^\text{th}$-order approximation to $u(\bt)$ for all targets $\bt \in
\overline{B_{h_{k}/2} (\bc_{j,k})}$, up to a precision given by the
quadrature error, discussed next.

\subsection{Controlling quadrature error
\label{sec:quaderror}}

Quadrature error in QBX is estimated in the second term in
Theorem~\ref{thm:qbx-orig}, i.e. it is mainly controlled by the
quantity $(h_k/r)^{2\qquad}$. In a way, this term measures the amount of
resolution supplied to numerically integrate singularities at a
distance $r$. Notionally, two mechanisms exist for controlling this
resolution--the panel length $h_k$ and the quadrature order $2\qquad$. Since
we have chosen to proportionally tie the distance of the singularity
$r$ to the panel length $h_k$, we do not expect the quadrature error
to change in response to mesh refinement, i.e. shrinking the panel
length $h_k$.  Thus the remaining instrument to control quadrature
error is the number of quadrature points $\qquad$ on each panel.
Below, we present an empirical procedure to determine a suitable
value for $\qquad$ for a given kernel to satisfy a user-specified
accuracy bound $\epsilon$.

Before we do so, we would like to highlight a family of situations that threaten
the accuracy of the evaluated layer potentials through an increase in
the quadrature error. Consider a source discretization consisting of
source panels with unequal panel lengths $h_\ell$ contributing to the
coefficients at a single expansion center $\bc_{j,k}$ belonging to a
`target' panel $\Gamma_k$.

If $\qquad$ in $(h_k/r)^{2\qquad}$ is chosen so as to only provide
sufficient quadrature resolution from panel $\Gamma_k$ to its
own centers $\bc_{j,k}$, then no variability in panel lengths $h_\ell$
can be tolerated without risking insufficient quadrature resolution being
available in some source panel/center combinations. Consider a panel
$\Gamma_\ell$ with length $h_\ell\approx 2 h_k$ adjacent to the panel
$\Gamma_k$ belonging to the current expansion center $\bc_{j,k}$.
While the center-to-source distance $r$ in $(h_\ell/r)^{2\qquad}$
typically does not obey a larger lower bound on $\Gamma_\ell$ than
$\Gamma_k$, the increase in the numerator leads to a worse error
estimate. Such a situation is common in adaptively refined meshes
which we wish to permit. To mitigate the impact of this effect,
we enforce the following condition:
\begin{cond}[Two-to-one length restriction between adjacent panels]
  \label{cond:two-to-one}
  If panel $\Gamma_k$ and panel $\Gamma_\ell$ are adjacent to each other, then
  \[
  h_{\ell} / h_{k} \in \left[1/2, 2 \right].
  \]
\end{cond}
We may then choose $\qquad$ so that $(2h_k/r)^{2\qquad}$ still
satisfies the required accuracy bound, and we therefore tolerate
2-to-1 panel sizing steps in adaptively refined meshes.

Unfortunately, the impact of differing panel sizes is not only felt in
adjacent panels. Rather, in the expansions at center $\bc_{j,k}$,
quadrature error
originating from, say, source panel $\Gamma_\ell$, is controlled by
$(h_\ell/d(\bc_{j,k}, \Gamma_\ell))^{2\qquad}$.
The left half of Figure~\ref{fig_nearby} illustrates a situation
in which inaccurate evaluation may occur.
If $(h_k/r)^{2\qquad}$
was chosen to satisfy accuracy constraints, then enforcing
\[
\frac{h_\ell}{d(\bc_{j,k}, \Gamma_\ell)}
\le
\frac{2h_k}{r}
\]
ensures that the error contribution of other source panels
$\Gamma_\ell$ will not exceed that of $\Gamma_k$.
Using our choice $r=h_k/2$, this leads to the following condition,
which we will algorithmically enforce.
\begin{cond}[Sufficient quadrature resolution from all source panels
  to all centers] \label{cond:quad-resolution}
  For each expansion center $\bc_{j,k}$, the distance
  $d(\Gamma_\ell, \bc_{j,k})$ from the center to all source panels $\Gamma_l$
  must at least be commensurate with the source panel's length $h_l$,
  or larger, to ensure adequate quadrature resolution:
  \[
    d(\bc_{j,k},\Gamma_\ell) \ge \frac{h_\ell}4
    \fourquad
    \text{for all expansion centers $\bc_{j,k}$ }
    \quad
    \text{for all source panels $\Gamma_l$}.
  \]
\end{cond}
We note that the task of enforcing Condition~\ref{cond:quad-resolution}
is somewhat challenging owing to its all-pairs (``all sources to all centers'')
nature. A fast (non-quadratic) algorithm for its enforcement will be presented in
Section~\ref{sec:quad-resolution-alg}.

In the setting of the Helmholtz equation, the
Helmholtz parameter $\omega^2$ (and the length scale of wave features
associated with it) represents a final aspect that may impact the
quadrature error. To avoid issues of this nature, we require that
source discretization panel lengths are bounded with respect to
this length scale, as expressed by the following condition.

\begin{cond}[Panel size bounded based on wavelength]
  \label{cond:wavelength-panel}
  The panel size is bounded with respect to the wavelength.
  \[
  \omega \cdot \left(\max_{k\in\{1,\dots,N\}} h_k \right) \leq 5.
  \]
\end{cond}

Returning to the choice of the number of Gauss quadrature points $\qquad$ per panel
to satisfy a user-specified relative accuracy bound
$\tilde C_{p,\qquad,\beta}(2h_k/r)^{2\qquad}<\epsilon$, we
note that Theorem~\ref{thm:qbx-orig} provides no explicit way to
estimate the constant $\tilde C_{p,\qquad,\beta}$. While work such as
\cite{af2016estimation} does provide explicit, non-asymptotic
formulas for the quadrature error, we note that these are, by
necessity, kernel-dependent. In the interest of generality, we
describe a numerical procedure, to be carried out once per kernel,
that finds a suitable number $\qquad$ for each given accuracy target
$\epsilon$.

\begin{figure}[t!]
\begin{center}
\includegraphics[width=.4\linewidth]{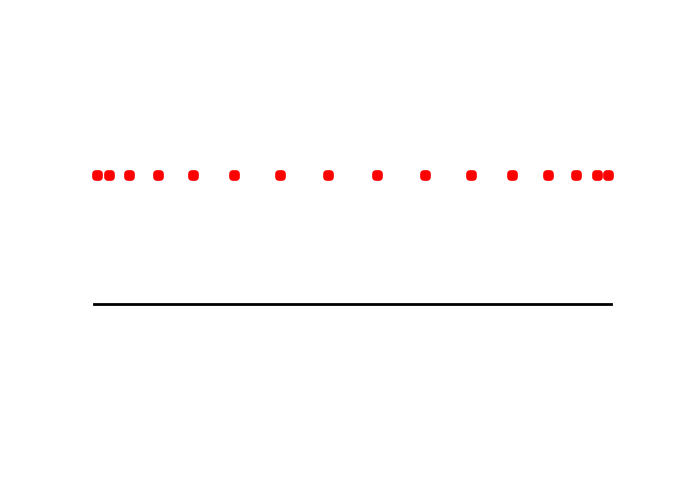}
\includegraphics[width=.4\linewidth]{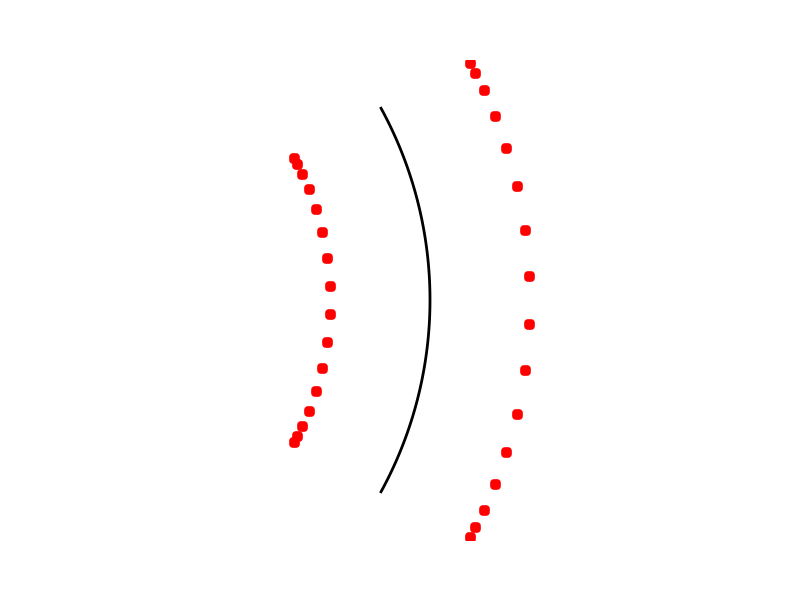}
\end{center}
\caption{Empirical estimation of the parameter $\qquad$.
Left: flat panel test geometry, Right: curved panel test geometry.
The curved panel is obtained as a unit length arc of a circle with
radius $1$.}
\label{geofig}
\end{figure}

We wish to convert the condition
\[\tilde C_{p,\qquad,\beta}\left(\frac{2h_k}r\right)^{2\qquad}<\epsilon\]
to one that can, at least approximately, be verified once for all panels
$\Gamma_k$. To this end, we first rescale so that $2h_k=1$, so that
an expansion center at distance $h_k/2$ would appear at a distance
$1/4$ from the unit-length panel. Next, we reduce to two `generic'
panel configurations, a straight panel and a curved panel
(cf.~Figure~\ref{geofig}). For all $\qtgt$ centers $\bc_j$ at a distance
of $1/4$ of the panel length, we ensure that the coefficient
integrals
\[
I_\ell(\bc_j) = \frac{i}{4} \int_{\Gamma} H_\ell^{(1)}( \omega
\vert\bc_{j} - \bx'\vert) \, e^{im\theta'} \, \sigma(\bx') \,
ds(\bx')
\]
are computed to within the specified accuracy $\epsilon$ by adjusting
$\qquad$.  As can be seen in the figure, for the straight panel, we
consider test centers on one side of the geometry (for symmetry
reasons), while we consider centers on both sides for the curved
one.

This test is accomplished through self-convergence, i.e.\ by
increasing $\qquad$ until the changes in the computed integrals in
response to resolution increase are below the specified threshold.
The tests are carried out with $\omega=5$
(cf. Condition~\ref{cond:wavelength-panel}) and with densities $\sigma
= P_n$, where $P_n$ is the $n^{th}$ degree Legendre polynomial for $n =
0,\ldots,\qtgt-1$.
Results of this experiment are summarized in Table~\ref{novertable},
for various choices of the accuracy parameter $\epsilon$.

\begin{table}[b!]
\caption{Source quadrature node count $\qquad$ as a
function of $\qtgt$ and $\epsilon$.}
\begin{center}
 \begin{tabular}{|c|c|c|c|c|}
 \hline
   & \multicolumn{4}{|c|}{$\epsilon$} \\ \hline
  $\qtgt$ & $10^{-3}$ & $10^{-6}$ & $10^{-9}$ & $10^{-12}$ \\ \hline
  2 & 8 & 16 & 24 & 32 \\ \hline
  4 & 12 & 24 & 32 & 40\\ \hline
  8 & 16 & 32 & 40 & 48\\ \hline
  16 & 32 & 48 & 64 & 64\\ \hline
 \end{tabular}
 \end{center}
 \label{novertable}
\end{table}

If the four conditions derived in this section are obeyed, then
the local expansion coefficients $\alpha_{\ell,m,n}$
for each QBX expansion centers $\bc_{m,n}$ calculated via
the oversampled discretization given by
\begin{equation}
\alpha_{\ell,m,n} =
\frac{i}{4} \sum_{k=1}^{N} \sum_{j=1}^{\qquad }
        \quadquant{w_{j,k}} \, H_{\ell}^{(1)} ( \omega \vert \bc_{m,n} -
                        \quadquant{\bs_{j,k}} \vert ) \, \quadquant{\sigma_{j,k}} \,
\end{equation}
for $\ell = -p,\ldots, p$ approximate the exact local
expansion coefficients of the single layer potential,
$\alpha_{\ell} (\bc_{m,n})$ defined in equation \eqref{eq:graf-coefficient},
to the prescribed tolerance $\epsilon$. In addition, we remark that
for locations \emph{not} covered by QBX expansion disks,
\emph{unmodified} Gaussian quadrature
\begin{equation}
u(\bt) = \mathcal{S}[ \sigma ] (\bt) \approx
\frac{i}{4} \sum_{k=1}^{N} \sum_{j=1}^{\qquad }
        \quadquant{w_{j,k}} \, H_{0}^{(1)} ( \omega \vert \bt -
                        \quadquant{\bs_{j,k}} \vert ) \, \quadquant{\sigma_{j,k}}
  \label{eq:gauss-layerpot}
\end{equation}
approximates the true potential $u(\bt)$ to the prescribed tolerance
$\epsilon$. We note that such targets $\bt$ satisfy
\[
d(\bt,\Gamma_{k}) > h_{k}/4 \fourquad \text{ for all $k$}.
\]
As a result, by evaluating QBX expansions where available and
using the `direct' Gauss-Legendre computation
\eqref{eq:gauss-layerpot} everywhere else, we obtain an
$\epsilon$-accurate approximation of the layer potential in
all of $\mathbb R^2$.


\subsection{Summary \label{sec:cons-summary}}

We have derived four sufficient conditions that, if satisfied,
guarantee sufficient accuracy in the evaluation of the layer
potential. They are summarized below:
\begin{description}
 \item [Condition~\ref{cond:no-sources-in-disk}] Expansion disk undisturbed by sources.
  \[d(\bc_{j,k},\Gamma_{\ell}) \geq h_{k} / 2\]
  for all panels numbers $\ell, k\in \{1,\dots,N\}$ and $j=1,\dots,\qtgt$.
 \item [Condition~\ref{cond:two-to-one}] Two-to-one length restriction between adjacent panels.
   If panel $k$ and panel $\ell$ are adjacent to each other, then
   \[h_{\ell} / h_{k} \in \left[1/2, 2 \right].\]

 \item [Condition~\ref{cond:quad-resolution}]
  Sufficient quadrature resolution from all source panels to all centers.
  \[d(\bc_{j,k}, \Gamma_{\ell}) \geq h_{\ell} / 4\]
  for all panels numbers $\ell, k\in \{1,\dots,N\}$ and $j=1,\dots,\qtgt$.
 \item [Condition~\ref{cond:wavelength-panel}] Panel size bounded based on wavelength.
  \[\omega \cdot \left(\max_{k\in\{1,\dots,N\}} h_k \right) \leq 5.\]
\end{description}
These conditions reduce the problem of accurate layer potential
evaluation to one of ensuring the validity of these four
geometric conditions. By inspection of each condition, it is clear
that if a given input discretization is found to be in violation of
any of the conditions, refining (i.e. splitting all offending panels
into two equal-length pieces) until all conditions are obeyed
is an effective remedy.

Conditions~\ref{cond:two-to-one} and \ref{cond:wavelength-panel} are
trivially checked by local computation, one panel at a time.
The other two conditions however are non-local in character, since
they concern pairs of (not necessarily related) sources and centers.
While they are trivial to verify by iterating over all such pairs,
doing so would negate the benefit of the fast algorithm being derived
in this contribution, since it would revert the run time of the
algorithm to scaling with the square of the input size, rather than
(nearly) linearly.

We finally note that the above family of conditions is technically
entirely independent of the FMM (or any other acceleration scheme). It simply
represents a family of conditions to ensure that Global QBX (in
which the expansions represent contributions from the entire source
geometry) to be accurate. It \emph{is} however natural to consider
the computational verification of these conditions in conjunction
with an FMM, since the tree data structures already available in this
setting can be put to excellent use.

\section{Area Queries in Quad-Trees}
\label{sec:area-query}%

Conditions~\ref{cond:no-sources-in-disk} and \ref{cond:quad-resolution}
of the previous section were found to involve non-local properties of
the input source discretization. We note that being able to
efficiently answer queries of the following type is useful
in the process of verifying both conditions:
\begin{quote}
  \textbf{Prototypical query (PQ):}
  Given a collection of $N$ points $\bx_{i}$, with an associated radius $r_{i}$,
  identify the collection of targets $\bt_{j}|_{j=1}^{M}$
  that are contained in $\cup_{i=1}^{N} B_{r_{i}}(\bx_{i})$.
\end{quote}
The goal of this section is to design an algorithm to complete precisely
this task in $\cO(N+M)$ CPU time.

Quad-tree data structures (such as those constructed by the Fast
Multipole Method \cite{greengard1998accelerating} base of the present
algorithm) are well-suited for such tasks.
Consider the following setup of an adaptive quad-tree.
Let $b_{0}$ be the smallest square centered at the origin which
contains all {\em particles}.
Particles may be comprised of any combination of sources, targets,
expansion centers, and centers of mass of all panels, depending
on the procedure.
We next introduce a hierarchy of meshes on the computational domain~$b_0$.
Mesh level~$0$ corresponds to~$b_0$ and mesh level~$\ell+1$ is
obtained by splitting boxes at level~$\ell$ into four quadrants,
denoted {\em children} of the {\em parent} box.
In order to allow for adaptivity, we allow for different levels
of refinement in different regions of $b_{0}$.
Empty boxes are pruned.
Let $\cB_{\ell}$ denote the set of non-empty boxes at level $\ell$.
We will refer to the result of this procedure as a quad-tree
on the computational domain $b_{0}$.

\begin{rem}
  In a quad-tree, a box $b$ is subdivided into $4$ equal parts if it
  contains more than some pre-specified number of particles,
  $n_{\text{max}}$.  It should be noted that not all categories of
  particles are considered for determining whether $b$ will be
  subdivided.  For example, in a quad tree containing all sources,
  targets, and expansion centers, $b$ can be partitioned if it
  contains more than $n_{\text{max}}$ sources, regardless of the
  number of targets, and expansion centers in $b$.  This subdivision
  criterion is decided depending on the procedure, and $n_\text{max}$
  is user-specified.
\end{rem}

To answer queries of type (PQ) given a source $\bx_{j}$ in
a leaf box $b$, it is \emph{not} sufficient to check all targets in box $b$
and the list of boxes adjacent to $b$ as the size of the box $b$ and
its adjacent boxes is independent of the extent of $\bx_{j}$.
It is also unsatisfactory to answer these queries by traversal from
the root, since placement of the query box near boundaries of the 
tree boxes will lead to expensive traversals of much of the tree
structure.

Instead, to accommodate these geometry queries,
it is sufficient to augment conventional adaptive
quad-tree data structures by adding the
capability to identify the minimal collection of leaf boxes which
completely contain a specified local neighborhood $B_{r_{j}} (\bx_{j})$.
Our proposed modification to the standard quad-tree consists of
an additional iteration structure available on the quad-tree, an
operation we term an \emph{area query}, defined as follows.
\begin{defin}[Area Query]
  \label{def:area-query}
  Given:
  \begin{itemize}
    \item A quad-tree partitioning of a square $b_0 \subseteq \mathbb R^2$, and
    \item a square $C_r(\bc)=\{\bx \in \mathbb R^2:
      \vert\bx-\bc\vert_\infty \le r\}$ with    $\bc\in b_0$,
  \end{itemize}
  an area query provides
  \begin{itemize}
    \item a list of all leaf (childless) boxes, $b_j$, in the quad-tree
      for which $b_j\cap C_r(\bc) \ne
    \emptyset$.
  \end{itemize}
\end{defin}
Note that while the point $\bc$ needs to lie in $b_0$, $C_r(\bc)
\subseteq b_0$ is not required.
A simple way to ensure that all relevant $\bc$ points are contained
in $b_0$ is to include them as points on which the quad-tree is built.

\subsection{Peers}

\begin{figure}
  \RawFloats
  \centering
  \begin{minipage}[t]{0.45\textwidth}
    \centering
    \begin{tikzpicture}[
        scale=0.4,
        fmmbox/.style={draw},
      ]
      \input{sample-tree.tex}

      \def\drawpeer#1{
        \fill [black!10] \boxpath{#1};
        \begin{pgfonlayer}{foreground}
          \node [fill=white,fill opacity=0.8,text opacity=1,font=\footnotesize] at
                                                (boxc#1) {$p$};
        \end{pgfonlayer}
      }
      \drawpeer{16}
      \drawpeer{17}
      \drawpeer{19}
      \drawpeer{7}
      \drawpeer{12}
      \drawpeer{24}
      \drawpeer{25}

      \fill [blue!40] \boxpath{18};
      \begin{pgfonlayer}{foreground}
        \node [fill=white,fill opacity=0.8,text opacity=1,font=\footnotesize] at
                                        (boxc18) {$b$};
      \end{pgfonlayer}

      \drawboxes

    \end{tikzpicture}
    \caption{
      The boxes marked $p$ are the peers of the box marked $b$.
    }
    \label{fig:peer-box}
  \end{minipage}\hfill
  \begin{minipage}[t]{0.45\textwidth}
    \centering
    \begin{tikzpicture}[
        scale=0.4,
        fmmbox/.style={draw},
      ]
      \input{sample-tree.tex}

      \def\drawpeer#1{
        \fill [black!10] \boxpath{#1};
      }
      \drawpeer{12}
      \drawpeer{14}
      \drawpeer{15}
      \drawpeer{7}
      \drawpeer{8}
      \drawpeer{10}
      \drawpeer{4}
      \drawpeer{13}

      \begin{pgfonlayer}{foreground}
        \node [fill=white,fill opacity=0.8,text opacity=1,font=\footnotesize]
          at (boxc13) {$b$};
      \end{pgfonlayer}

      \drawboxes

      \def\drawaq#1{
        \fill [black!30]
          ($ (boxl#1) + (0.1, 0.1) $)
          rectangle
          ($ (boxh#1) - (0.1, 0.1) $);
      }
      \drawaq{4}
      \drawaq{22}
      \drawaq{31}
      \drawaq{39}
      \drawaq{38}
      \drawaq{18}
      \drawaq{33}
      \drawaq{35}
      \drawaq{32}
      \drawaq{34}
      \drawaq{27}
      \drawaq{24}
      \drawaq{26}

      \coordinate (pt) at ($ (boxc13) + (-0.14*\boxsizeXIII,0.4*\boxsizeXIII) $);
      \def\rad{\boxsizeXIII * 0.45}
      \fill (pt) circle (2pt);
      \draw [thick] ($ (pt) - (\rad,\rad) $) rectangle ($ (pt) + (\rad,\rad) $) ;

      \draw[ultra thick,dashed,blue!80] \boxpath{13};
    \end{tikzpicture}
    \caption{%
      An area query associated with center $\bc$ denoted
      as a thick dot. The extent of the area query is
      shown by a thick outline.
      The guiding box $b$ associated with this area query is
      marked by a dashed blue line.
      The peers of box $b$
      are shown in light gray, and the results of the area query are shown with
      a dark gray inset.
    }
    \label{fig:guiding-box}
  \end{minipage}
\end{figure}

While the area query is conceptually easy to understand, care must be
taken that its actual implementation is efficient.
To this end, we first define the notion of  a \emph{peer box}.
\begin{defin}[Peer Box]
  Given a box $b_j$ in a quad-tree, $b_k$ is a \emph{peer box} of $b_j$ if
  it is
  \begin{enumerate}
    \item adjacent to $b_j$,
    \item of at least the same size as $b_j$
      (i.e.\ at the same or a coarser level), and
    \item no child of $b_k$ satisfies the above two criteria.
  \end{enumerate}
  \label{def:peer-box}
\end{defin}
Recall from the original adaptive
FMM~\cite{carrier_fast_1988} that a \emph{colleague} is an
adjacent box at the same level. A \emph{peer box} is similar, however
a box at a coarser level may be included in the set of peer boxes if a
colleague at the same level fails to exist.
Figure~\ref{fig:peer-box} shows an example.
The box $b_j$ itself is included in its list of peers.
Any box (in two dimensions) has at most nine peers.

\subsection{Performing Area Queries}
\label{sec:area-query-algorithm}

To process a single area query for the region $C_r(\bc)$, as in
Definition~\ref{def:area-query}, we first determine
the {\em guiding box} associated with the region
$C_r(\bc)$.
Referring to Figure ~\ref{fig:guiding-box},
the guiding box $b$ for the query, is the smallest
box $b$ which contains $\bc$ such that
$C_{r} (\bc)$ is completely contained in the peers
of $b$.
The algorithm below identifies all leaf boxes
resulting from the area query associated with
$C_r(\bc)$.
\qline
\begin{center}
  {\bf Algorithm for area queries}
\end{center}
{\small
\begin{tabbing}
......\=......\=......\=......\=      \kill
\noindent {\bf Comment} [Determine guiding box $b$ of query]\\
\> Set $b$ to be the root box $b_0$.\\
\> \textbf{do}\\
\>\> \textbf{if} $\vert b \vert/2 < r \le \vert b\vert$, \\
\>\>\> \textbf{break} out of the loop,\\
\>\> \textbf{else if} $b$ has a child containing $\bc$,\\
\>\>\> find the child $\tilde b$ of $b$ containing $\bc$, set $b =
\tilde b$,\\
\>\> \textbf{else}\\
\>\>\> \textbf{break} out of the loop,\\
\>\> \textbf{end if}\\
\> \textbf{end do}\\
\noindent {\bf Comment} [Enumerate and check leaf descendants of
peers of $b$ ]\\
\> \textbf{do} $b_j\in\{b_k:b_k\text{ is a peer of $b$}\}$\\
\>\> \textbf{do} $b_\ell\in\{b_k:b_k\text{ is a childless
descendant of $b_j$}\}$\\
\>\>\> \textbf{if} $b_\ell\cap C_r(\bc) \ne \emptyset$,\\
\>\>\>\> \textbf{include} leaf box $b_\ell$ in the resulting set for the area query.\\
\>\>\> \textbf{end if}\\
\>\> \textbf{end do}\\
\> \textbf{end do}
\end{tabbing}
}
\qline
In the above algorithm and in the following, by $\vert b_j \vert$ we mean the \emph{radius}
of the box $b_j$, i.e.\ the axis-aligned distance from the center of the
box to its edge. For example, a unit box $b$ has size $\vert b \vert = 1/2$.
Practical implementations can be designed to take advantage of
concurrency and locality of data access by processing a large number of
area queries in batch form.

The area query mechanism improves on a simple generalization of
the quad tree where {\em macro}-sources have a specified
interaction extent.
Area queries ensure that the
length-scale of the macro-source is commensurate with the box
structure being used for geometric look-up.
The three main advantages of area queries are:
(a) it is trivially parallelizable, (b) although presented here in two
dimensions, it generalizes directly to three-dimensional geometries,
and (c) it
is competitive in operation count with other algorithmic options that we have
explored for detecting regions of validity for QBX expansions.

\subsection{Complexity and Correctness}

The two following lemmas are straightforward to prove using arguments
based on the structure of standard quad-trees.

\begin{lem}
  The area query algorithm of Section~\ref{sec:area-query-algorithm}
  performs $O(L+n_{\text{leaf}})$ work, where $L$ is the number of levels in
  the tree, and $n_{\text{leaf}}$ is the number of leaf boxes being returned by
  the query.
\end{lem}

\noindent
This complexity estimate is readily apparent from the algorithm
 above, as there are at most $L$ iterations of the loop locating the
{\em guiding box} $b$, the number of peers of a box is bounded by a
constant, and the number of descendants of $b$ exceeds the number of boxes
returned by at most a constant factor.

\begin{lem}
  A box is returned by the  area query algorithm of
  Section~\ref{sec:area-query-algorithm} if and only if it
  satisfies Definition~\ref{def:area-query}.
\end{lem}
If a box $b_n$ is returned from the query, it necessarily is a leaf
box overlapping $C_r(\bc)$. Conversely, let $b_n$ be a leaf box overlapping
$C_r(\bc)$.   Thus $|\bc_n - \bc|_\infty \le r + |b_n|$. Let $b_i$ be
the area query's guiding box. Then
$|\bc_i-\bc|_\infty \le |b_i|$. Note that $|b_n|\le |b_i|$ and $r\le |b_i|$.
Combining these facts yields
\[
  |\bc_i-\bc_n|
  \le |\bc_i-\bc| + |\bc-\bc_n|
  \le |b_i| + r+ |b_n| \le 3|b_i|.
\]
Therefore, $b_i$'s peers cover at least the area given by
$\{\bx\subset \mathbb R^2:\vert \bx-\bc_i \vert_\infty \le 3|b_i|\}$,
so $b_n$ must be a child of one of the peers of $b_i$ and thus was
examined and returned by the area query.

\section{Triggering Source Refinement}
\label{sec:source-refine}

Using the area query algorithm of the previous section, we next
describe algorithms to verify the conditions of
Section~\ref{sec:cons-summary}. In addition to merely detecting
violations, the algorithms further attribute the detected issues to one
or more panels $\Gamma_k$ which may then be refined. Upon refinement,
the conditions are rechecked, and, if necessary, additional rounds of refinement
and checking are performed, until all conditions are satisfied.

On refinement, a flagged panel
is subdivided  into two panels of equal arclength, determined via the
Legendre expansion describing the panel's parametrization.
As discussed earlier, Conditions~\ref{cond:two-to-one} 
and~\ref{cond:wavelength-panel} are easily checked one panel at a time.
In this section, we outline algorithms to flag panels based on
Conditions~\ref{cond:no-sources-in-disk} and~\ref{cond:quad-resolution}.

\subsection{An Algorithm to Verify Condition~\ref{cond:no-sources-in-disk}}
\label{sec:no-source-in-disk-alg}

For notational convenience, we refer to expansion centers $\bc_{j}$
and sources $\stgt_j$ with a single index.
For determining panels which violate Condition~\ref{cond:no-sources-in-disk},
if $\bmm_k$ is the center of mass of panel $\Gamma_k$, then let
$r_{k}$ be the smallest radius such that the entire panel $\Gamma_k$
is contained in an $\ell^\infty$ `disk' $Z_k$ of radius $r_k$
centered at $\bmm_k$, i.e.
\begin{equation}
\label{eq:tkdef}
\Gamma_{k}
\subseteq \{\bx: \vert \bx - \bmm_{k} \vert_{\infty} \leq r_{k}\} = Z_k.
\end{equation}

Next, let $S(b)$ denote the list of panels overlapping a leaf box $b$,
so that $\Gamma_{k} \in S(b)$ if $\Gamma_{k} \cap b \neq \emptyset$.
It is straightforward to construct this list using
an area query on $\bmm_{k}$ and its associated bounding square $Z_{k}$.

Let $\bc_{k}$ be the expansion center associated with node
$\stgt_k$ on panel $m$.
Let $r_{\bc_{k}} = h_{m}/2$, and
let $C_{k} = \{ \bx \in \mathbb{R}^{2} :
\vert\bx - \bc_{k}\vert_\infty \le r_{\bc_{k}}\}$ denote the area query
search domain associated with $\bc_{k}$.
For each leaf box $b$ returned by the area query, we loop
over all panels in $S(b)$ and flag panels which violate
Condition~\ref{cond:no-sources-in-disk}.
The algorithm below identifies all such panels.

\qline
\begin{center}
  {\bf Algorithm for triggering refinement based on
  Condition~\ref{cond:no-sources-in-disk}}
\end{center}
{\small
\begin{tabbing}
......\=......\=......\=......\=......\=      \kill
 {\bf Comment} [Choose main parameters]\\
\> Create a quad-tree  on the computational domain containing
        all expansion centers and centers of mass of all panels.\\
\> Choose the maximum number $n_\text{max}$ of particles in a childless box.\\
\> Subdivide a box $b$ if it contains more than $n_{\text{max}}$
expansion centers.
\end {tabbing}
\begin{center}
          {\bf Stage 1.}
\end{center}
\noindent {\bf Comment} [Refine the computational cell into a hierarchy
of meshes for sorting expansion centers, and centers of mass of all
panels.]
\begin{tabbing}
......\=......\=......\=......\=......\=      \kill
\> {\bf do $\ell = 0,1,2,\ldots$}     \\
\> \> {\bf do } $b_k\in \cB_{\ell}$ \\
\> \> \> {\bf if} $b_k$ contains more than $n_\text{max}$
particles {\bf then}  \\
\> \> \> \>  subdivide $b_k$ into four boxes, ignore (prune) the
empty boxes formed.     \\
\> \> \>{\bf end if}   \\
\> \> {\bf end do}  \\
\>  {\bf end do} \\
\\
 {\bf Comment} [Let $n_\text{box}$ be the total number of boxes.] \\
 {\bf Comment} [Let $M(b)$ and $E(b)$ denote the
list of panel centers of mass and expansion centers in box $b$
respectively.]
\end {tabbing}
\begin{center}
{\bf Stage 2.}
\end{center}
\begin{tabbing}
......\=......\=......\=......\=......\=......\=......\=      \kill
 {\bf Comment} [Loop over all panel centers of mass in every leaf box.
Using area queries, identify the list of panels $S(b)$ \\
relevant for each leaf box $b$.]\\
\\
\> {\bf do} $j = 1,2,\ldots, n_\text{box}$     \\
\> \> {\bf if} $b_j$ is childless {\bf then}  \\
\> \> \>  {\bf do $\bmm_{k} \in M( b_{j})$} \\
\> \> \> \> Perform an area query for the region $Z_{k}$
to identify the list of leaf boxes $A_k$. \\
\> \> \> \> {\bf do $b_\ell \in A_k$} \\
\> \> \> \> \> Append panel $\Gamma_{k}$ to $S(b_\ell)$ \\
\> \> \> \> {\bf end do} \\
\> \> \>{\bf end do}   \\
\> \> {\bf end if}  \\
\>  {\bf end do}
\end{tabbing}
\begin{center}
          {\bf Stage 3.}
\end{center}
\begin{tabbing}
......\=......\=......\=......\=......\=......\=......\=      \kill
 {\bf Comment} [Loop over all expansion centers in every leaf box.
Using area queries, loop over all relevant \\
 panels to find expansion
centers which violate condition~\ref{cond:no-sources-in-disk} and
flag the corresponding
panel.]\\
\\
\> {\bf do} $j = 1,2,\ldots, n_\text{box}$     \\
\> \> {\bf if} $b_j$ is childless {\bf then}  \\
\> \> \>  {\bf do $\bc_{k} \in E( b_{j})$} \\
\> \> \> \> Perform an area query for the region $C_{k}$
to identify the list of leaf boxes $A_k$. \\
\> \> \> \> {\bf do $b_\ell \in A_k$} \\
\> \> \> \> \> {\bf do $\Gamma_{m} \in S(b_\ell)$} \\
\> \> \> \> \> \> {\bf if} $d(\bc_{k},\Gamma_{m}) \leq h_{n}/2$ where
                        $\bc_{k}$ is an expansion center associated with
                                  a source on panel $n$, \\
\> \> \> \> \> \>       {\bf and} $m\neq n$, {\bf then} Flag panel $n$. \\
\> \> \> \> \> {\bf end do} \\
\> \> \> \> {\bf end do} \\
\> \> \>{\bf end do}   \\
\> \> {\bf end if}  \\
\>  {\bf end do}
\end {tabbing}
}
\qline

\subsection{An Algorithm to Verify Condition~\ref{cond:quad-resolution}}%
\label{sec:quad-resolution-alg}%
To determine panels which violate Condition~\ref{cond:quad-resolution},
we carry out the following algorithm.
If $\bmm_\ell$ is (again) the center of mass of panel $\Gamma_\ell$, then let
$\tilde r_{\ell}$ be the smallest radius such that
\begin{equation} \label{eq:tldef}
T_\ell = \{\bx: d(\bx,\Gamma_{\ell}) \leq
h_{\ell}/4 \} \subseteq \{\bx: \vert \bx - \bmm_{\ell} \vert_{\infty}
\leq \tilde r_{\ell}\} = \tilde C_\ell.
\end{equation}
Note that $T_\ell$ is a tubular neighborhood around panel $\Gamma_\ell$,
and
$\tilde C_\ell$ is a bounding square centered at the panel's center of mass
containing $T_\ell$.

\begin{figure}
  \def\quadrefinebasepic{
      \draw
        (0.26, 3.62) coordinate (n00)
        .. controls +(0:0.7cm) and +(180:0.7cm) ..
        (2.00, 4.15) coordinate (n01)
        .. controls +(0:0.7cm) and +(180:0.7cm) ..
        coordinate [pos=0.4] (coarsetgt)
        coordinate [pos=0.5] (coarsecmass)
        coordinate [pos=0.1] (inacc1)
        coordinate [pos=0.7] (inacc2)
        (3.61, 3.92) coordinate (n02)
        .. controls +(0:0.7cm) and +(180:0.7cm) ..
        (4.93, 4.18) coordinate (n03)
        node [anchor=west] {$\Gamma_{\text{coarse}}$};

      \draw [yshift=0.6cm]
        (0.92, 1.90+0.15) coordinate (n10)
        .. controls +(0:0.3cm) and +(180:0.3cm) ..
        (1.87, 2.03+0.15) coordinate (n11)
        .. controls +(0:0.3cm) and +(180:0.3cm) ..
        coordinate [pos=0.3] (finetgt)
        coordinate [pos=0.5] (finecmass)
        (2.78, 1.90+0.15) coordinate (n12)
        .. controls +(0:0.3cm) and +(180:0.3cm) ..
        (3.68, 2.0+0.15) coordinate (n13)
        .. controls +(0:0.3cm) and +(180:0.3cm) ..
        (4.52, 1.94+0.15) coordinate (n14)
        node [anchor=west] {$\Gamma_{\text{fine}}$};

    \foreach \i in {0,...,3}
      \fill (n0\i) circle (1pt);

    \foreach \i in {0,...,4}
      \fill (n1\i) circle (1pt);

    \draw (finetgt) circle (1pt);

    \def\finerad{0.5}
    \def\coarserad{0.9}

    \path (finetgt) ++(80:\finerad) coordinate (finectr);
  }
  \begin{tikzpicture}[scale=1.4]
    \quadrefinebasepic

    \draw [dashed] (finetgt) -- (finectr);

    \fill (finectr) circle (1pt);

    \draw [red] (finectr) circle (\finerad);

    \draw [thick,red,shorten >=0.2cm,->]
      (inacc1)
      .. controls +(250:0.5cm) and +(110:0.5cm) ..
      (finectr);
    \draw [thick,red,shorten >=0.2cm,->]
      (inacc2)
      .. controls +(250:0.5cm) and +(50:0.5cm) ..
      (finectr)
      node [pos=0.5,anchor=north west] {inaccurate!};

  \end{tikzpicture}
  \hfill
  \begin{tikzpicture}[scale=1.4]
    \quadrefinebasepic

    \draw (coarsecmass) circle (1pt);

    \draw [dashed] (finetgt) -- ++(80:\finerad) coordinate (finectr);

    \fill (finectr) circle (1pt) node [anchor=west] {$\bc_{j,k}$};

    \begin{pgfonlayer}{background}
      \fill[black!20] (n01) -- ++(0,-1)
        .. controls +(0:0.7cm) and +(180:0.7cm) ..
        ($ (n02) + (0,-1) $) -- (n02)
        -- ++(0,1)
        .. controls +(180:0.7cm) and +(0:0.7cm) ..
        coordinate [pos=0.5] (tklabel)
        ($ (n01) + (0,1) $);
    \end{pgfonlayer}

    \draw [|<->|] (n01) -- +(0,1)
      node [pos=0.5,anchor=east] {$h_\ell/4$};
    \def\aqrad{1.12}
    \draw [|<->|] (coarsecmass) -- +(\aqrad,0)
      node [pos=0.5,anchor=south] {$\tilde r_\ell$};

    \node [anchor=west,text=blue] at ($ (coarsecmass) + (\aqrad, \aqrad*0.5) $)
      { $\tilde C_\ell$};

    \node [below=1mm of tklabel] {$T_\ell$};
    \node [above=0mm of coarsecmass] {$\Gamma_\ell$};
    \node [below=0mm of coarsecmass] {$\bmm_\ell$};

    \draw [blue]
      ($ (coarsecmass) + (-\aqrad,-\aqrad) $)
      rectangle
      ($ (coarsecmass) + (\aqrad,\aqrad) $);


  \end{tikzpicture}

  \caption{Refinement triggered by nearby (but not adjacent)
    geometry. On the left, the potential due to
    $\Gamma_{\text{coarse}}$ at the QBX expansion center near $\Gamma_{\text{fine}}$
    would not be computed accurately. The local panel size
    on $\Gamma_{\text{coarse}}$ is too large to obey the correct error
    estimate.
    On the right, an area query in region $C_k$, in which the
    panel length on $\Gamma_{\text{coarse}}$ triggers refinement of
    $\Gamma_{\text{coarse}}$.
  }
  \label{fig_nearby}
\end{figure}
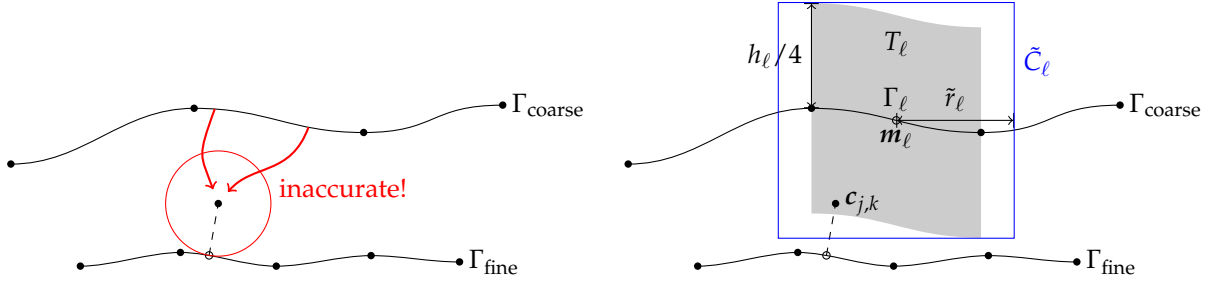

If an expansion center on panel $k$ lies in $T_\ell$
defined in equation \eqref{eq:tldef}, with $k \neq
\ell$, and panel $k$ is not adjacent to panel $\ell$,
then panel $\ell$ will be flagged.
$\tilde C_{\ell}$ will be the area query search domain associated with $\Gamma_{\ell}$
to search for centers $\bc_{j,k}$ which may receive insufficiently
resolved quadrature contributions from $\Gamma_\ell$.
The algorithm below identifies all such panels.
The left-hand side of Figure~\ref{fig_nearby} illustrates the general situation
being detected, the right-hand side clarifies the notation.

\qline
\begin{center}
  {\bf Algorithm for triggering refinement based on
  Condition~\ref{cond:quad-resolution}}
\end{center}
{\small
\noindent {\bf Comment} [Choose main parameters]
\begin{tabbing}
......\=......\=......\=......\=      \kill
\> Create a quad-tree  on the computational domain containing
        all expansion centers and centers of mass of all panels.\\
\> Choose the maximum number $n_\text{max}$ of particles in a childless box.\\
\> Subdivide a box $b$ if it contains more than $n_{\text{max}}$
expansion centers.
\end {tabbing}

\begin{center}
          {\bf Stage 1.}
\end{center}

\noindent {\bf Comment} [Refine the computational cell into a hierarchy
of meshes for sorting expansion centers, and centers of mass of all
panels.]
\begin{tabbing}
......\=......\=......\=......\=      \kill
\> {\bf do $\ell = 0,1,2,\ldots$}     \\
\> \> {\bf do $b_k \in \cB_\ell$ } \\
\> \> \> {\bf if} $b_k$ contains more than $n_\text{max}$ particles {\bf then}  \\
\> \> \> \>  subdivide $b_k$ into four boxes, ignore (prune) the empty
             boxes formed.     \\
\> \> \>{\bf end if}   \\
\> \> {\bf end do}  \\
\>  {\bf end do}
\end {tabbing}

\noindent {\bf Comment} [Let $n_\text{box}$ be the total number of boxes.] \\
\noindent {\bf Comment} [Let $M(b)$ and $E(b)$ denote the
list of panel centers of mass and expansion centers in box $b$, respectively.]
\begin{center}
          {\bf Stage 2.}
\end{center}

\noindent {\bf Comment} [Loop over centers of mass of all panels
in every leaf box.
Using area queries, loop over all relevant expansion centers to find expansion
centers which violate Condtion \ref{cond:quad-resolution} and flag the corresponding panel.]

\begin{tabbing}
......\=......\=......\=......\=......\=......\=      \kill
\> {\bf do} $j = 1,2,\ldots, n_\text{box}$     \\
\> \> {\bf if} $b_j$ is childless {\bf then}  \\
\> \> \>  {\bf do $\bmm_{k} \in M(b_j)$} \\
\> \> \> \> Perform an area query for the region $\tilde C_{k}$
to identify the list of leaf boxes $A_k$ (cf. Figure~\ref{fig_nearby}, right panel)\\
\> \> \> \> {\bf do $b_\ell \in A_{k}$} \\
\> \> \> \> \> {\bf do $\bc_{m} \in E(b_\ell)$} \\
\> \> \> \> \> \>
{\bf if} $d\left(\bc_{m},\Gamma_{k}\right) \leq h_{k}/4$, where $\bc_{m}$ is an
expansion center associated with a source node on panel $n$, \\
\> \> \> \> \> \>
with $n \neq k$ and panel $n$ is not adjacent to panel $k$ {\bf then} flag
panel $n$. \\
\> \> \> \> \> {\bf end do} \\
\> \> \> \> {\bf end do} \\
\> \> \>{\bf end do}   \\
\> \> {\bf end if}  \\
\>  {\bf end do}
\end {tabbing}
}
\qline


\section{Evaluating Layer Potentials in the Volume} 
\label{sec:volume-eval}

QBX-type expansions of the potential can be used not only to
accurately evaluate on-surface, but also near-surface potentials
\cite{barnett_evaluation_2012}. It is thus expedient to
present an algorithm that similarly unifies the treatment of potential
evaluation at \emph{any} target, regardless of whether it
is located on, near, or far from the source geometry.
This section presents an efficient geometric algorithm based on area
queries that enables the subsequent fast algorithm to provide
this capability.

\subsection{Marking Targets for Evaluation by QBX}
\label{target-flagging}

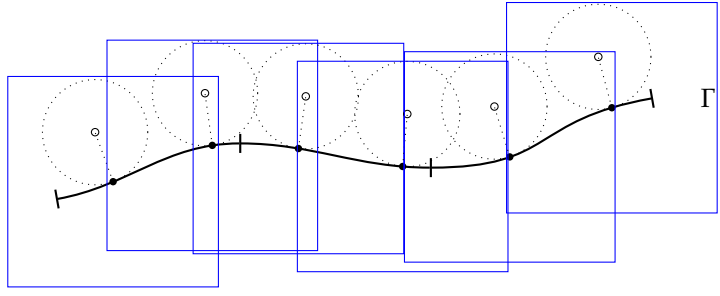
\begin{figure}
  \def\marklegendre#1{
      \foreach \i / \ln in {
        2/0.29707742,
        5/0.87076559,
      }
      { coordinate [pos=\ln] (#1\i) }
    }

  \centering
  \begin{tikzpicture}[scale=1.4]
    \coordinate (a) at (0.26, 3.62);
    \coordinate (b) at (2.00, 4.15);
    \coordinate (c) at (3.81, 3.92);
    \coordinate (d) at (5.93, 4.58);

    \draw [|-,thick]
      (a) .. controls +(10:0.7cm) and +(180:0.7cm) .. (b)
      \marklegendre{a};

    \draw [|-,thick]
      (b) .. controls +(0:0.7cm) and +(180:0.7cm) .. (c)
      \marklegendre{b};

    \draw [|-|,thick]
      (c) .. controls +(0:1.2cm) and +(190:1.2cm) .. (d)
      node [anchor=west,xshift=5mm] {$\Gamma$}
      \marklegendre{c};

    \foreach \intv in {a,b,c}
      \foreach \nd in {2,5}
        \fill (\intv\nd) circle (1pt);

    \def\rad{0.5}
    \foreach \nd/\angle in {
      a2/110,
      a5/98,
      b2/82,
      b5/85,
      c2/107,
      c5/105
      }
    {
      \draw [dotted] (\nd) -- ++(\angle:\rad cm) coordinate (ctra\nd);

      \draw (ctra\nd) circle (1pt);

      \draw [dotted] (ctra\nd) circle (\rad);
      \draw [blue] (\nd) +(-2*\rad,-2*\rad) rectangle +(2*\rad,2*\rad);
    }

  \end{tikzpicture}
  \caption{Area queries (blue) being carried out to locate
  off-geometry/volume targets (not shown)
  that are close enough to the source geometry to require potential
  evaluation through a QBX expansion. The location of QBX centers and
  the regions where their expansions are valid are shown by dotted
  lines.
  One area query is being carried out per source point
  on the source geometry. To simplify the image, only two target
  points (with associated area queries) are shown per geometry panel.
  }
  \label{fig:tunnel-sources}
\end{figure}

The objective of this section is to describe an
algorithm to determine, for each target point $\bt_\ell$,
whether a `conventional' evaluation of the potential based on
Gaussian quadrature and the (point-based) FMM is sufficient, or whether
evaluation of the potential through a QBX expansion is required
(and if so, which
of the many available centers should be chosen).
As discussed in Section~\ref{sec:quaderror}, the potential at
targets which satisfy $d (\bx,\Gamma_{k}) \leq h_{k}/4$ needs to be
evaluated through a QBX expansion.
We will refer to this tubular domain {\em close} to the boundary as
$\Gamma_\text{near}$ given by:
\begin{equation}
  \begin{aligned}
    \Gamma_\text{near} &= \cup_{k=1}^{N}
    \left\{\bx: \, d(\bx,\Gamma_{k}) \leq h_{k}/4 \right\} \\
&= \cup_{k=1}^N T_k.
\end{aligned}
\end{equation}
For efficiency, the algorithm proceeds in three stages:
\begin{enumerate}
  \item Set up `tunnel' area queries around source geometry to mark region
    in which potentially inaccurate layer potential evaluation
    can occur.
    Figure~\ref{fig:tunnel-sources} illustrates the area query being performed.
  \item Tag leaf boxes covered by area queries with
    `endangering' source geometry
  \item Based on leaf box containment and tagging information
    from previous step, decide (a) whether the potential at a
    center needs to be evaluated through QBX and (b) whether
    a corresponding center is available.
\end{enumerate}

\qline
\begin{center}
  {\bf Associate volume/surface targets with QBX expansions}
\end{center}
\begin{tabbingalg}
\noindent {\bf Comment} [Set up source area query]\\
\> \textbf{do} $k=1,\dots,N$ (panels)\\
\>\> \textbf{do} $j=1,\dots,\qtgt$ (sources on panel $k$)\\
\>\>\> Find bounding squares $S_{j,k,\pm}$ containing the disks about centers
$\bc_{j,k}$ with radii $\frac{h_k}4$\\
\>\>\> Find bounding square $S_{j,k}$ centered at $\squad_{j,k}$ so that $S_{j,k}\supseteq S_{j,k,+}\cup S_{j,k,-}$\\
\>\> \textbf{end do}\\
\> \textbf{end do}\\
\end{tabbingalg}
\begin{tabbingalg}
\noindent {\bf Comment} [Mark boxes with nearby sources]\\
\>Perform an area query over the squares $S_{j,k}$\\
\> \textbf{do} in each of the resulting leaf boxes $b_i$\\
\>\> Add $\squad_{j,k}$ to the set $D_i$\\
\> \textbf{end do}\\
\end{tabbingalg}
\begin{tabbingalg}
\noindent {\bf Comment} [Decide target association]\\
\> \textbf{do} for each target $\ell=1,2,\ldots, n_{t}$\\
\>\> Find leaf box $b_i$ containing $\bt_\ell$ \\
\>\> \textbf{if} $\bt_\ell \in \Gamma_{\text{near}}$ \\
\>\>\> Mark $\bt_\ell$ as requiring evaluation by QBX\\
\>\> \textbf{end if}\\
\>\> {\bf Comment} [Locate closest eligible QBX center]\\
\>\> $(j', k')=\operatorname{argmin}_{(j,k)}
  \{ |\bt_\ell - \bc_{j,k}| : \squad_{j,k}\in D_i,
 |\bt_\ell - \bc_{j',k'}| \le \frac{h_{k'}}2 (1+\epsilon_{\text{assoc}}) \}$\\
\>\> \textbf{if} $(j',k')$ were found \textbf{then}\\
\>\>\> \textbf{mark} $t_\ell$ to use the QBX expansion around $\bc_{j',k'}$\\
\>\> \textbf{else}\\
\>\>\> \textbf{if} $\bt_\ell$ was marked as requiring QBX\\
\>\>\>\> \textbf{fail} target association\\
\>\>\> \textbf{else}\\
\>\>\>\> \textbf{mark} $t_\ell$ for evaluation without QBX\\
\>\>\> \textbf{end if}\\
\>\> \textbf{end if}\\
\> \textbf{end do}
\end{tabbingalg}
\qline
The tolerance $\epsilon_{\text{assoc}}$ in the algorithm is used to
ensure that a target point $\bt_\ell$ located on the source geometry
$\Gamma$ is properly associated with QBX centers on the source
geometry in the presence of inexact arithmetic.

In addition, it should be remarked that the expansion disks do not fully
cover the area immediately surrounding the source geometry. Some gaps
remain. Algorithmically, we could approach this issue via refinement or by
adding additional centers.
Empirically, increasing $\epsilon_{\text{assoc}}$ to cover
any targets in this region leads to no decrease in accuracy. This is
mathematically at least plausible since the expanded layer potentials
are very smooth even in the immediate neighborhood of the source
geometry, and hence the terms of the QBX expansion grow
slowly.

It is straightforward to see that the area query algorithm finds a
superset of the targets located in $\Gamma_{\text{near}}$. To
determine whether a target $\bt_\ell$ is actually located within
$\Gamma_{\text{near}}$, one may employ Newton's method to find the
closest point on the panels $\Gamma_k$, or as an approximation,
one may use
\[
d(\bt_\ell, \Gamma_k)\approx
\min_j |\bt_\ell - \stgt_{j,k}|_2.
\]
If accurate evaluation of the layer potential on \emph{both}
sides of the geometry is desired, it is straightforward to augment
the presented scheme with a {\em side-preference} mechanism that restricts
eligible QBX centers to ones on a predetermined side of the geometry.
This is particularly important for points located \emph{on} the
source geometry, since for these targets it is impossible to determine
which of the two limits is desired by geometric location alone.

\section{A Fast Algorithm for QBX} 
\label{sec:fmm}

Given a tolerance, $\epsilon$, and the number original of Gauss-Legendre
 points per panel,
$\qtgt$, we determine $\qquad$ from Table~\ref{novertable}.
As a pre-processing step, we interpolate the discretized
geometry and density from the density grid to the source grid
As in~\eqref{eq_sigdisc},
we then approximate the layer potential $u = \mathcal S [\sigma]$
using the sum
\begin{equation}
 u(\bx) \approx
 \frac{i}{4} \sum_{k=1}^{N} \sum_{j=1}^{\qquad}
 \quadquant{w_{j,k}} \, H^{(1)}_{0}( \omega
\lvert \bx - \squad_{j,k} \rvert) \,
\sigquad_{j,k}.
\end{equation}
For notational convenience, we rewrite the above sum as
\begin{equation}
 u(\bx) \approx \frac{i}{4}
 \sum_{j=1}^{n_{s}} \quadquant{w}_j \, H^{(1)}_0 (\omega \vert \bx - \squad_j
 \vert) \, \sigquad_j
 \label{eq:QBXFMMpot}
\end{equation}
where $n_{s} = N\qquad$.
We also approximate the $J$-expansion coefficients $\alpha_{\ell,j,k}$ at the
expansion center $\bc_{j,k}$ using the same source-grid
Gauss-Legendre quadrature rule:
\begin{equation}
 \alpha_{\ell,j,k} \approx \frac{i}{4} \sum_{n=1}^{N}
 \sum_{m=1}^{\qquad} \quadquant{w}_{j,k} \,
H^{(1)}_\ell (\omega\vert \bc_{j,k} - \squad_{m,n}\vert) \,
e^{i\ell\theta'} \, \sigquad_{m,n} ,
\end{equation}
for $\ell = -p,\ldots,p$, $j=1,\ldots,q$, and
$k=1,\ldots,N$.
Again for notational convenience, we rewrite the above expression as
\begin{equation}
 \alpha_{\ell,j} \approx \frac{i}{4} \sum_{j=1}^{n_{s}} \quadquant{w}_{j}
 \, H^{(1)}_\ell (\omega\lvert
 \bc_j - \squad_{j} \rvert) \, \sigquad_{j},
 \label{eq:localexpcoeffs}
\end{equation}
for $j=1,\ldots,Nq$.
In a minor abuse of notation,
we will interchange $\approx$ and $=$ when discussing discrete sums.
The task at hand is to accelerate the computation of:
\begin{enumerate}
 \item the potential $u(\bt_j)$
 defined in equation~\eqref{eq:QBXFMMpot} at the target locations $\bt_{j}$
which are not flagged to be in $\Gamma_{\text{near}}$ and
 \item the $J$-expansion coefficients $\alpha_{\ell,j}$ for each
   expansion center $\bc_j$, defined in equation ~\eqref{eq:localexpcoeffs}.
\end{enumerate}
The FMM has traditionally been used to accelerate
the computation of the potential $u(\bt)$ defined
in equation~\eqref{eq:QBXFMMpot}.
Roughly speaking, the algorithm heirarchically compresses
the {\em far-field} interactions which are numerically
low-rank.
We describe below a conceptually and algorithmically simple modification to the original
FMM algorithm, to accelerate the far-field interactions
in the computation of the local expansion coefficients
$\alpha_{l,j}$ defined in equation~\eqref{eq:localexpcoeffs}.

In a Helmholtz FMM,
based on $\epsilon$, we determine $\pfmm \approx
\log(\epsilon)$, the multipole expansion order
for the $H$-expansions and $J$-expansions of the FMM.
In practice, this parameter can vary depending on which level the
translations inside the FMM are being processed.
Let $b_{0}$, the computational domain, be the smallest square
centered at the origin which contains all expansion centers, sources,
and targets.
Assume that $b_{0}$ is partitioned using a
quad-tree, and that for any box $b$ in the tree, let $F(b)$ denote
the far-field of the box $b$.
The far-field of a box $b$ is the collection of boxes
which are {\em well-separated} from the box $b$ at the length-scale
of the size of $b$.
By $\psi_b$ we denote the $J$-expansion for box $b$:
\begin{equation}
\psi_{b}(\bx) = \sum_{\ell=-\pfmm}^{\pfmm} \gamma_\ell
\,  J_\ell (\omega\lvert \bx - \bmm_{b}\rvert) \,
e^{-i \ell \theta_{\bx,\bmm_{b}} } \, ,
\end{equation}
where it is assumed that $\bx \in b$, $\bmm_b$ is the center of box $b$, and
in polar form,
$\bx - \bmm_{b} = ( \rho, \theta_{\bx,\bmm_{b}})$.
The expansion $\psi_b$
is an $\epsilon$-approximation to the potential due to to all sources that
are in \mbox{$F(b)$}:
\begin{equation}
\left\lvert  \psi_{b} (\bx) -
\frac{i}{4}  \sum_{\squad_{j} \in F(b)} \quadquant{w}_j \,
 H^{(1)}_0 (\omega\vert\bx - \squad_{j}\vert) \, \sigquad_{j}
\right\rvert
= \mathcal O(\epsilon).
\end{equation}
For a particular expansion center $\bc_{j}$ contained in $b$, by using the
standard $J$-expansion to $J$-expansion (local-to-local) translation
 operator
we can obtain a $J$-expansion of order $p \leq \pfmm$ about $\bc_{j}$
given by $\tilde{\psi}_{c_{j}}$ as
\begin{equation}
 \tilde{\psi}_{c_{j}}(\bx) = \sum_{\ell=-p}^{p}
 \tilde{\gamma}_{\ell,j} \,
 J_\ell (\omega\vert \bx - \bc_{j} \vert) \, e^{-i \ell \theta_{\bx,\bc_{j}} } \, ,
\end{equation}
where the local polar coordinates are given by
$\bx - \bc_{j} = ( \rho, \theta_{\bx,\bc_{j}})$.

Using Graf's Addition Theorem for $H^{(1)}_0$, we see that
$\tilde{\gamma}_{\ell,j}$ corresponds to the contribution to
$\alpha_{\ell,j}$ from sources $\squad_{k} \in F(b)$:
\begin{align}
 \tilde{\psi}_{c_{j}} &= \sum_{\ell=-p}^{p} \tilde{\gamma}_{l,j} \,
 J_\ell (\omega |\bx - \bc_{j} | ) \,
        e^{-i \ell \theta_{\bx,\bc_{j}} } \, , \\
 &= \frac{i}{4} \sum_{ \squad_{k} \in F(b)}
 H^{(1)}_0  (\omega |\bx - \squad_{k} | ) \,
         \sigquad_{k} + \mathcal O(\epsilon )  , \\
 &= \sum_{\ell=-p}^{p} \left( \frac{i}{4} \sum_{ \squad_{k}\in F(b)}
 H^{(1)}_\ell (\omega |\squad_{k} - \bc_{j} | ) \,
 e^{i \ell \theta_{\squad_{k},\bc_{j}} } \, \sigquad_{k} \right)
 J_\ell (\omega |\bx - \bc_{j} |) \, e^{-i \ell \theta_{\bx,\bc_{j}} }
   + \mathcal O(\epsilon ) .
\end{align}
Therefore, we have that
\begin{equation}
\tilde{\gamma}_{l,j} = \frac{i}{4}  \sum_{\squad_{k} \in F(b)}
 H^{(1)}_\ell (\omega |\squad_{k} -\bc_{j} | ) \,
         e^{i \ell \theta_{\squad_{k}, \bc_{j}} } \,
         \sigquad_{k} + \mathcal O(\epsilon)  .
\end{equation}

From the above discussion, it is easy to see that we can accelerate
far-field computation of
the potentials at un-flagged target locations as well as
the $J$-expansion
coefficients at expansion centers using small modifications of a
standard FMM.
For those already familiar with FMMs, we will describe briefly the
modifications
required to the standard (`point') FMM. For a detailed description of the algorithm,
we refer the reader to Section~\ref{sec:complete-fmm}.

In order to compute values of the potential at un-flagged targets,
we do not need to make any modifications to the standard FMM.
To compute the $J$-expansion coefficients at the expansion centers,
we need the following four additional steps.
Using the standard notation for FMM interaction lists for a box $b$,
$U(b)$ and $W(b)$,
(see Section~\ref{sec:complete-fmm} for a detailed
definition of these lists),
for an expansion center $\bc$ in a leaf box $b$
of the tree hierarchy:
\begin{enumerate}
 \item Form the $J$-expansion due to all sources
   $\squad_j \in U(b)$,
 \item Form the $J$-expansion by translating the $H$-expansion of all
         boxes  $b'\in W(b)$ to account for all sources
         $\squad_j \in W(b)$,
 \item Translate the $J$-expansion of the box $b$ to a $J$-expansion
 at $\bc$ to  account for all sources $\squad_j \in F(b)$,
\item Add the above three $J$-expansions together.
\end{enumerate}
To evaluate the potential at targets in $\Gamma_{\text{near}}$ which were
flagged in Section~\ref{target-flagging}, we use the $J$-expansion
of the corresponding  expansion center that was computed above.


\begin{figure}[t]
\begin{floatrow}
\ffigbox{
 \includegraphics[width=7cm]{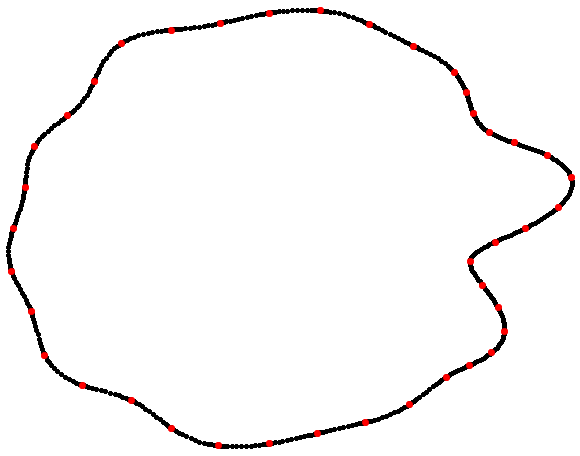}
}
{
        \caption[Test geometry for determining $p_{add}$]
{Test geometry for determining $p_{add}$}
\label{fig:test-padd}
}
\capbtabbox
{
 \begin{tabular}{|@{\hskip 8pt}c@{\hskip 8pt}|
         @{\hskip 10pt}c@{\hskip 10pt}
         |@{\hskip 10pt}c@{\hskip 10pt}|
                @{\hskip 8pt}c@{\hskip 8pt}|
                @{\hskip 8pt}c@{\hskip 8pt}|}
  \hline
  $p$ & 2 & 4 & 6 & 8 \\ \hline
  $p_{add}$ & 5 & 5 & 15 & 20 \\ \hline
 \end{tabular}
}
{
        \caption{$p_{add}$ as a function of $p$ }
 \label{paddtable}
}
\end{floatrow}
\end{figure}

\subsection{Maintaining Expansion Accuracy}

The order of the multipole
expansion in the FMM, $\pfmm$, is dependent on the tree-level
$\ell$
and frequency $\omega$.
Standard estimates are available for choosing this parameter for
$H$-expansions and $J$-expansions (see, for
example~\cite{greengard1998accelerating}).
Briefly, the expansion order $\pfmm$ for an outgoing expansion is chosen
to evaluate the sum
\begin{equation}
\sum_{j=1}^{N} H_{0}^{(1)} (\omega \vert \bx-\squad_{j} \vert) \sigma_{j}
\label{eq:pointfmmsum}
\end{equation}
for $\bx$ in the far-field of the sources, with error less than $\varepsilon$.
Suppose all the sources
$\squad_{j}$ are contained in a box centered at the origin with
$\vert b \vert = R$, where $R$ is the size of the box at level
$\ell$.
Using Graf's addition theorem, the outgoing expansion
associated with the box $b$ is given by
\begin{equation}
\sum_{k=-\infty}^{\infty} a_{k} H_{k}^{(1)} (\omega \rho) e^{ik\theta}\, ,
\end{equation}
where
\begin{equation}
a_{k} = \sum_{j=1}^{N} J_{k} (\omega \rho_{j}) e^{-ik\theta_{j}} \sigma_{j} \,.
\end{equation}
Here $(\rho,\theta)$ and $(\rho_{j},\theta_{j})$
are the polar coordinates of $\bx$ and $\squad_{j}$ respectively.
The targets in the far-field of the box $b$ are
are separated from box $b$ by at least one box length, i.e. $\rho\geq 3R$.
In order to compute the sum~\ref{eq:pointfmmsum} with precision $\varepsilon$,
the outgoing expansion is truncated at $\pfmml$
\begin{equation}
\sum_{k=-\infty}^{\infty} a_{k} H_{k}^{(1)}(\omega\rho) e^{ik\theta}
= \sum_{k=-\pfmml}^{\pfmml} a_{k} H_{k}^{(1)}(\omega\rho) e^{ik\theta}
+ \mathcal{O}(\epsilon) \, ,
\end{equation}
if
\begin{equation}
\max_{\substack{\rho\geq 3R \\ \rho_{j} \leq \sqrt{2}R \\ \vert n\vert > \pfmml}}
\vert H_{n} (\omega \rho) J_{n} (\omega \rho_{j}) \vert \leq \epsilon \, .
\end{equation}

In the QBX framework, the expansion order $\pqbx$ must be chosen to compute
the sums
\begin{equation}
\sum_{j=1}^{N} H_{m}^{(1)} (\omega \vert \bx -\squad_{j} \vert) e^{im\theta}
\sigma_{j}
\, , \label{eq:qbxfmmsum}
\end{equation}
for all $\vert m \vert \leq p$, and $\bx$ in the far-field
of the sources, with error less than $\epsilon$.
Using Graf's addition theorem, the outgoing expansion corresponding to
the sum in Equation~\eqref{eq:qbxfmmsum} is given by
\begin{equation}
\sum_{k=-\infty}^{\infty} a_{k} H_{k+m}^{(1)} (\omega \rho)
e^{i(k+m)\theta}\, , \label{eq:outexpqbx}
\end{equation}
where
\begin{equation}
a_{k} = \sum_{j=1}^{N} J_{k} (\omega \rho_{j}) e^{-ik\theta_{j}} \sigma_{j} \,.
\end{equation}
Thus, outgoing expansion~\eqref{eq:outexpqbx} can be truncated at $\pqbxl$ if
\begin{equation}
\max_{\substack{\rho\geq 3R \\ \rho_{j} \leq \sqrt{2}R \\ \vert n\vert > \pqbxl
}}
\vert H_{m+n} (\omega \rho) J_{n} (\omega \rho_{j}) \vert \leq \epsilon \, ,
\end{equation}
for all $\vert m \vert \leq p$.
A similar analysis can be done for the incoming expansions as well.

While the above explanation provides an intuition for the need
of larger outgoing and incoming expansions for the FMM-accelerated
QBX as compared to the standard FMM, a detailed analysis for estimating $\pqbxl$ is
fairly involved.
In the evaluation of the discretized layer potential $\cS[\sigma]$,
the far-field of the sum~\eqref{eq:qbxfmmsum}
is scaled by $J_{m} (\omega h/2)$,
where $h$ is a characteristic arc-length of a panel
in the discretization of the boundary.
Moreover, the size of the smallest box in the quad-tree
data structure is also intricately tied to $h$.
Thus, we set $\pqbxl = \pfmml + \padd$ and
determine $\padd$ numerically, as a function of
$p$ and $\epsilon$, by testing
the accuracy of the translated $J$-expansions at targets close to the
boundary for $1000$ random geometries.
The Helmholtz parameter $\omega$ was set to $5$ for these numerical
experiments.
The boundaries of the random test geometries were
described by the following parametrization:
\begin{align}
 x_{1}(\theta) &= r(\theta) \cos\theta  ,\\
x_{2}(\theta) &= r(\theta) \sin\theta  ,
\end{align}
where
\begin{equation}
 r(\theta) = 5.0 + \sum_{j=1}^{12} \delta_{j} \, \sin j \theta  \, ,
\end{equation}
where $\delta_{j}$ are uniformly distributed in $\left[-0.2, 0.2 \right]$ and
$\theta \in [0,2\pi)$ (see, Figure~\ref{fig:test-padd}, for example).

\begin{figure}[t]
  \centering
  \begin{subfigure}[t]{.45\linewidth}
    \centering
    \includegraphics[width=.9\linewidth]{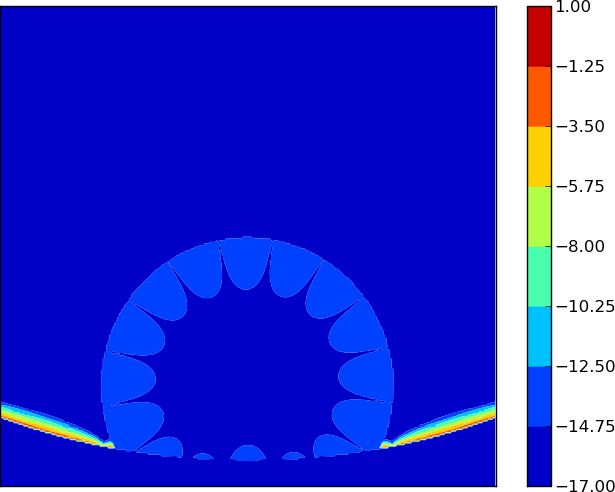}
    \caption{Directly computing $J$-expansion coefficients at the expansion
    center.}
  \end{subfigure}
  \fourquad
  \begin{subfigure}[t]{.45\linewidth}
    \centering
    \includegraphics[width=.9\linewidth]{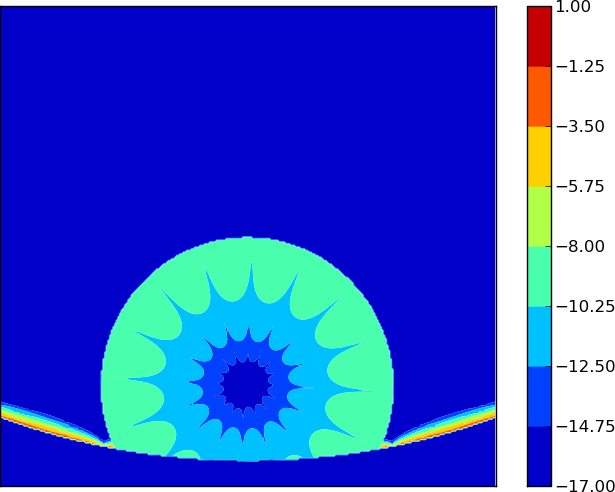}
    \caption{Using FMM-accelerated QBX with $p_{add} = 0$.}
  \end{subfigure}\\
  \begin{subfigure}[t]{.45\linewidth}
    \centering
    \includegraphics[width=.9\linewidth]{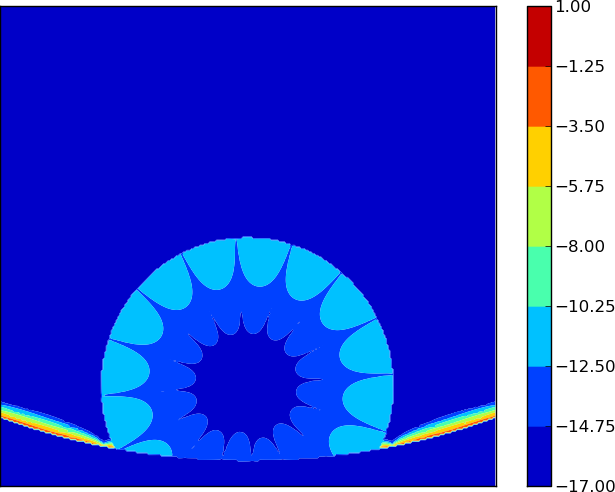}
    \caption{Using FMM-accelerated QBX with $p_{add} = p$.}
  \end{subfigure}
  \fourquad
  \begin{subfigure}[t]{.45\linewidth}
    \centering
    \includegraphics[width=.9\linewidth]{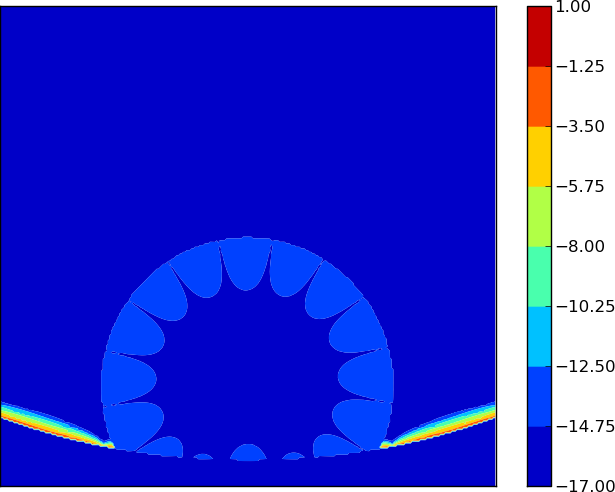}
    \caption{Using FMM-accelerated QBX with $p_{add} = 20$.}
  \end{subfigure}
  \caption{Error at targets in
    $\Gamma_{\text{near}}$ for
    $p = 8$.}
  \label{fig:padd}
\end{figure}

In particular,
 analogous to the subsequent
numerical experiments in Section~\ref{sec:acc-comp},
given $\epsilon$ and $p$, we compute
the relative error in Green's identity for a known Helmholtz potential
at a collection of targets close to the boundary.
The associated layer potentials are evaluated using the FMM-accelerated
QBX algorithm described in \ref{sec:complete-fmm}, and
the geometry is sufficiently refined and over-sampled to ensure that
the relative error in Green's identity using a direct calculation is
less than $\epsilon$.
The order of the multipole expansion in the FMM is set to
$\pfmml +p'$, where we vary the parameter $p'$.
The resulting $\padd$ is the minimum $p'$ for which the relative error in
Green's identity computed using the FMM-accelerated QBX algorithm is
less than $\epsilon$.
A contour plot of the error in evaluating the layer potential close
to the boundary, wherein the local expansion coefficients are computed
using a direct computation, and an FMM-accelerated QBX scheme for
different values of $\padd$ is shown in Figure~\ref{fig:padd}.
The experiments indicated that $\padd$ was independent of the prescribed
precision $\epsilon$, and merely a function of $p$.
The results are summarized in Table~\ref{paddtable}

\subsection{Complete Statement of the Fast Algorithm}
\label{sec:complete-fmm}
To complement the previous discussion, and for mathematical and
algorithmic completeness, we provide in this section a complete
statement of the fast algorithm, including its `point' and `layer
potential' parts, but excluding the geometry preprocessing described
in detail earlier, with the differences to its prior `point-only' version
highlighted.
Let $b_{0}$ be the smallest square centered at the origin which contains all
sources, targets and expansion centers.
Let $|b_{0}|$ be half the length of the side of the root-box square (its
`radius').
Following the procedure described in Section~\ref{sec:area-query-algorithm},
construct a quad-tree on $b_{0}$.

In a minor abuse of notation, for any subset $A$ of the computational
box, let $\squad \in A$, $\bt \in A$, and $\bc \in A$
denote the set of sources, targets, and expansion centers contained
in $A$, respectively.
Let $\ell_{max}$ be the highest level of refinement at any point.

A box $b$ is a {\em parent box} if it has been subdivided into one or
more boxes.

A {\em child box} is a non-empty box resulting from the subdivision of
a parent box.

{\em Colleagues} are adjacent boxes at the same level including the self-box.
A given box has at most nine colleagues.

A {\em leaf box} is a childless box.

Let $\bmm_{b}$ denote the coordinates of the center of box $b$.
Boxes $b$ and $b'$ at level $\ell$ are {\em well-separated} from each other if
\begin{equation}
 | \bmm_{b} - \bmm_{b'} | \geq 2 \cdot (2|b_0|)
\cdot 2^{-\ell} .
\end{equation}

The {\em $U$-list} of a box $b$, denoted by $U(b)$, is empty if $b$ is
a parent box.
If $b$ is a leaf box, $U(b)$ is the set of all leaf boxes that
are adjacent to $b$.

The {\em $V$-list} of a box $b$, denoted by $V(b)$, consists of
all the children
of the colleagues of the parent of $b$ that are well-separated from $b$.

The far-field of a box $b$ will be denoted by $F(b) = b_0 \setminus
(U(b) \cup W(b))$.

The {\em $W$-list} of a box $b$, denoted by $W(b)$, is empty if $b$ is a
parent box.
If $b$ is a leaf box, $W(b)$ consists of all the descendants of the
colleagues of $b$ whose parents are adjacent to $b$, but who are not
adjacent to $b$ themselves.
Note that a box $b' \in W(b)$ is separated from $b$ by a distance
equal to the length of the side of $b'$.

The {\em $X$-list} of a box $b$, denoted by $X(b)$, is formed by
all boxes $b'$ such that $b\in W(b')$.
Note that all boxes in the $X$-list are childless and larger than $b$.

Let $\phi_{b}$ denote the $\pqbxl$-term $H$-expansion about the center
of $b$ of the potential created by all sources in $b$.

Let $\psi_{b}$ denote the $\pqbxl$-term $J$-expansion about the center of
box $b$ of the potential created by all sources in the far-field of $b$,
that is $\squad \in F(b)$. The value
$\psi_{b}(\bt)$ is the result of evaluating the $J$-expansion at target~$\bt$.

Let $\Upsilon_{b}$ denote the $J$-expansion about the center
of $b$ representing the potential
due to all sources $\squad \in V(b)$.

Let $\Delta_{b}$ denote the $J$-expansion about the center of $b$
representing the potential due to all sources
$\squad \in X(b)$.

Let $\alpha_{b}(\bt )$ denote the potential at
$\bt \in b$ due to all sources $\squad \in U(b)$.

Let $\beta_{b}(\bt )$ denote the potential at
$\bt \in b$ due to all sources $\squad \in W(b)$.

Let $\tilde{\alpha}_{j,\ell} $ denote denote the $\ell^\text{th}$ $J$-expansion
coefficient at $\bc_{j}\in b$ due to all sources
$\squad \in U(b)$.

Let $\tilde{\beta}_{j,\ell}$ denote the $\ell^\text{th}$ $J$-expansion coefficient
at $\bc_{j} \in b$ due to all sources $\squad \in W(b)$.

Let $\tilde{\gamma}_{j,\ell}$ denote the $\ell^\text{th}$ $J$-expansion coefficient
at $\bc_{j} \in b$ due to all sources $\squad \in
b_{0} \setminus ( U(b) \cup W(b) )$.\\

\qline
\begin{center}
  {\bf FMM-accelerated QBX}
\end{center}
{\small
\noindent {\bf Comment} [Choose main parameters]
\begin{tabbing}
......\=......\=......\=......\=      \kill

\> Given $\epsilon$, using standard multipole estimates
\cite{rokhlin1990rapid}, set number of terms in expansions at level
$\ell$ to $\pfmml$. \\
\> Depending on $p$, use Table~\ref{paddtable} to determine $\padd$ and
set $\pqbxl = \pfmml + \padd$. \\
\> Create a quad-tree  on the computational domain containing
        all sources, targets, and expansion centers.\\
\> Choose the maximum number $n_\text{max}$ of particles in a childless box.\\
\> Subdivide a box $b$ if the sum of the number of sources and
targets in $b$ is greater than $n_{\text{max}}$.\\
\end {tabbing}

\noindent {\bf Comment} [Refine the computational cell into a hierarchy
of meshes.]
\begin{center}
          {\bf Stage 1.}
\end{center}

\begin{tabbing}
......\=......\=......\=......\=      \kill
\> {\bf do} $\ell = 0,1,2,\cdots$     \\
\> \> {\bf do} $b_j \in \cB_\ell$ \\
\> \> \> {\bf if} $b_j$ contains more than $n_\text{max}$
particles {\bf then}  \\
\> \> \> \>  subdivide $b_j$ into four boxes, ignore the empty boxes
             formed, add the non-empty boxes formed to $\cB_{\ell+1}$. \\
\> \> \>{\bf end if}   \\
\> \> {\bf end do}  \\
\>  {\bf end do}
\end {tabbing}
\noindent {\bf Comment} [Let $n_\text{box}$ denote the total number of boxes.]
\begin{center}
          {\bf Stage 2.}
\end{center}

\noindent {\bf Comment} [For every box $b$ at every level $\ell$,
        form a multipole
 expansion representing  the potential outside $b$ due to all
the particles contained in $b$.]\\




\noindent {\bf Comment} [For each childless box $b$, combine all
charges inside $b$ to obtain the $H$-expansion about
the center of $b$.]
\begin{tabbing}
......\=......\=......\=......\=......\=      \kill
\> {\bf do} $j=1, n_\text{box}$ \\
\> \>  {\bf if} $b_j$ is a childless box {\bf then} \\
\> \> \> form a $\pqbxl$-term $H$-expansion, $\phi_{b_j}$
representing the potential  outside $b_j$ due \\
\> \> \> to all charges located in $b_{j}$.       \\
\> {\bf end do}
\end {tabbing}




\noindent {\bf Comment } [For each parent box $b$, obtain the multipole
        expansion $\phi_{b}$ by translating the
$H$-expansions centered $\bmm_{b'}$ to an $H$-expansion centered
at $\bmm_{b}$, where
$b'$ is a child of $b$. Add the resulting expansions together.]

\begin{tabbing}
......\=......\=......\=......\=      \kill
\> {\bf do} $\ell=\ell_\text{max}-1,\ldots,1$  \\
\> \>   {\bf do} $b_j \in \cB_\ell$       \\
\> \> \> {\bf if} $b_j$ is a parent box {\bf then}                     \\
\> \> \> \> For each child of $b_j$, shift the center of the
$H$-expansion to $b_j$'s center. \\
\> \> \> \> Add the resulting expansions together to obtain
the expansion $\phi_{b_j}$.  \\
\> \> \> {\bf end if}                                            \\
\> \> {\bf end do}                                          \\
\> {\bf end do}
\end {tabbing}

\begin{center}
          {\bf Stage 3.}
\end{center}

\noindent {\bf Comment} [For all particles in each childless box $b$, compute
the interactions with all sources $\squad \in U(b)$  directly.]
\begin{tabbing}
......\=......\=......\=......\=      \kill
\> {\bf do} $j=1, n_\text{box}$           \\
\> \> {\bf if} $b_j$ is childless {\bf then}     \\
\> \> \>  For each target $\bt$ in $b_{j}$, compute the sum $\alpha_{b}(\bt)$
          of the interactions
          between $\bt$ and all sources
                                        $\squad \in U(b_{j})$. \\
\> \> \> {\bf New:} For each expansion center $\bc_{k}$ in $b_{j}$, compute
         the $J$-expansion \\
\> \> \> coefficients, $\tilde{\alpha}_{k,\ell}$ for $\ell=-p,\ldots p$,
         due to all sources $\squad \in U(b_{j})$. \\
\> \> {\bf end if}                                     \\
\> {\bf end do}
\end {tabbing}

\begin{center}
          {\bf Stage 4.}
\end{center}

\noindent {\bf Comment} [For each box $b$, convert
the $H$-expansions of all boxes in $V(b)$ into  $J$-expansions about
the center\\ of box $b$.]
\begin{tabbing}
......\=......\=......\=......\=      \kill
\> {\bf do} $j=1, n_\text{box}$ \\
\> \> {\bf do $b_k \in V(b_j)$}      \\
\> \> \> Convert $H$-expansion $\phi_{b_k}$ centered at $\bmm_{b_{k}}$
into a $J$-expansion centered at $\bmm_{b_{j}}$. \\
\> \> \> Add the resulting expansions to obtain $\Upsilon_{b_j} $.  \\
\> \> {\bf end do}                                           \\
\> {\bf end do}
\end {tabbing}

\begin{center}
          {\bf Stage 5.}
\end{center}

\noindent {\bf  Comment}
[For each childless box $b$,
 evaluate the $H$-expansions of all boxes in  $W(b)$
 at every particle position in $b$.]
\begin{tabbing}
......\=......\=......\=......\=      \kill
\> {\bf do} $j=1, n_\text{box}$                                   \\
\> \> {\bf if} $b_j$ is childless {\bf then}               \\
\> \> \> Evaluate the  $H$-expansion $\phi_{b_k}$ of each box
         $b_k \in W(b_j)$ to obtain $\beta_{b_j}(\bt)$
                                 for every target $\bt$ in $b_j$.\\
\> \> \> {\bf New:} Convert the $H$-expansion $\phi_{b_{k}}$ of each
        box $b_{k} \in W(b_j)$ to obtain the $J$-expansion\\
\> \> \> coefficients $\tilde{\beta}_{m,\ell}$, for $\ell=-p \ldots p$,
for each expansion center $\bc_{m} \in b_{j}$. \\
\> \> {\bf end if}                                                \\
\> {\bf end do}
\end {tabbing}

\begin{center}
          {\bf Stage 6.}
\end{center}

\noindent {\bf  Comment}
[For each  box $b$, form local
 expansions about the center $\bmm_b$ representing the potential
due to all\\ sources $\squad \in X(b)$.]
\begin{tabbing}
......\=......\=......\=......\=      \kill
\> {\bf do} $j=1, n_\text{box}$ \\
\> \>  Convert the potential of all sources
$\squad \in X(b_j)$ into a $J$-expansion about the
center of $b$.   \\
\> {\bf end do}
\end {tabbing}

\begin{center}
          {\bf Stage 7.}
\end{center}

\noindent {\bf Comment} [Shift the centers of $J$-expansions
of parent boxes to the centers of their children.]
\begin{tabbing}
......\=......\=......\=......\=      \kill
\> {\bf do} $\ell=1,\ell_{max}-1$  \\
\> \> {\bf do $b_j \in \cB_\ell$ }   \\
\> \> \> {\bf if } $b_j$ is a parent box {\bf then}   \\
\> \> \> \> Shift the center of expansion $\Upsilon_{b_j}$ to the center
of each of $b_{j}$'s children $b_{k}$.\\
\> \> \> \>  Add the resulting expansion to $\Upsilon_{b_k}$. \\
\> \> \> {\bf end if}  \\
\> \> {\bf end do}  \\
\> {\bf end do}
\end {tabbing}

\begin{center}
          {\bf Stage 8.}
\end{center}

\noindent {\bf Comment} [For each childless box $b$, obtain $\psi_b$
as the sum of local  expansions $\Upsilon_{b}$ and $\Delta_{b}$.
For each target $\bt$ in a childless box $b$, evaluate $\psi_{b}(\bt)$
and obtain the potential at $\bt$ by adding $\psi_{b}(\bt)$,
 $ \alpha_{b}(\bt)$ and   $\beta_{b}(\bt)$  together.]
\begin{tabbing}
......\=......\=......\=......\=      \kill
\> {\bf do} $j=1,n_\text{box}$ \\
\> \> {\bf if} $b_j$ is childless {\bf then}      \\
\> \> \> Compute $\psi_{b_j}=\Upsilon_{b_j}+\Delta_{b_j}$. \\
\> \> \> For each target $\bt$ in $b_j$, evaluate $\psi_{b_j}(\bt)$.\\
\> \> \> Add  $\psi_{b_j}(\bt)$, $\alpha_{b_j}(\bt)$ and $\beta_{b_j}(\bt)$
         to obtain the potential at $\bt$. \\
\> \> \> {\bf New:} For each expansion center $\bc_k$ in box $b_{j}$,
translate  $\psi_{b_j}$ to
compute $J$-expansion\\
\> \> \> coefficients $\tilde{\gamma}_{k,\ell}$, $\ell=-p,\ldots,p$
about $\bc_{k}$. \\
\> \> \> Add $\tilde{\alpha}_{k,\ell}$, $\tilde{\beta}_{k,\ell}$ and $\tilde{\gamma}_{k,\ell}$ to
obtain the $J$-expansion coefficient $\alpha_{k,\ell}$ at expansion center $\bc_k$. \\
\> \> {\bf end if} \\
\> {\bf end do}
\end {tabbing}
}
\qline
\section{Numerical Results} 
\label{sec:results}
In the following subsections
we illustrate the performance of FMM-accelerated QBX,
both in terms of accuracy and speed.
We demonstrate accuracy in evaluating layer potentials at targets
both on the boundary and in the volume by verifying Green's identity
using known solutions to the Helmholtz equation.
We also show the linear-time complexity of the algorithm described
in Section~\ref{sec:complete-fmm}, and compare its computational performance
with timings for the underlying standard point-FMM.
Finally, we use the Global QBX algorithm to solve a large
multi-scale scattering
problem with over 100,000 unknowns.
All the experiments in this section are performed using a single core
on a Dell laptop with a 2.2 GHz \mbox{Intel Core~i5-5200U}
processor and $8$~GB of RAM.
The \texttt{gfortran} compiler, version~4.9.3, was used.

\subsection{Preliminaries \label{sec-num-prel}}
\begin{figure}[t]
  \centering
  \begin{subfigure}[t]{.45\linewidth}
    \centering
    \includegraphics[width=.95\linewidth]{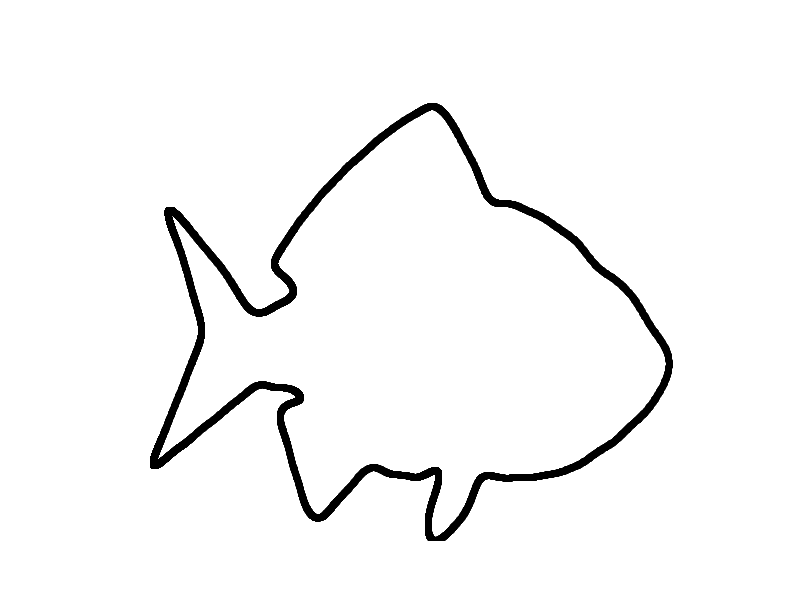}
    \caption{A sample fish-like geometry, analytically parametrized
      as a Fourier series.}
    \label{fig:sample_geom0}
  \end{subfigure}
  \fourquad
  \begin{subfigure}[t]{.45\linewidth}
    \centering
    \includegraphics[width=.95\linewidth]{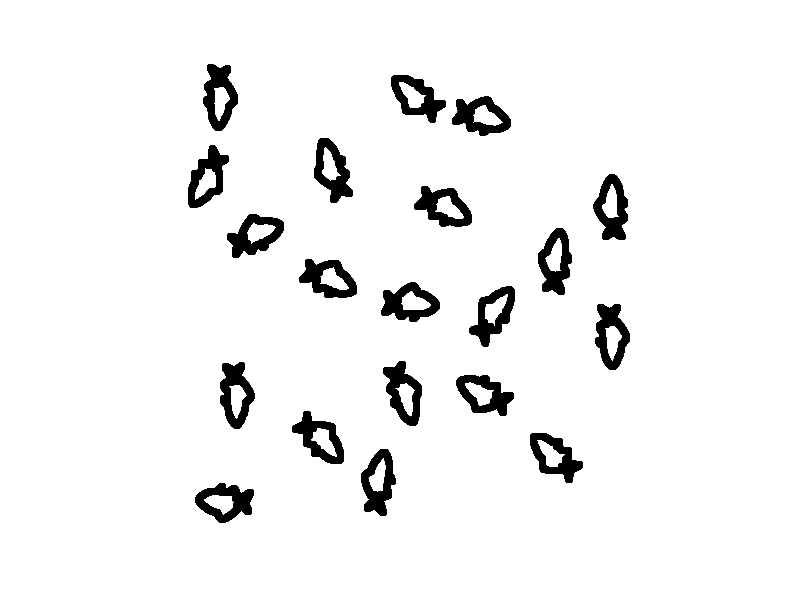}
    \caption{Sample geometry of rotated and translated fish
      for performance analysis.}
    \label{fig:sample_geom}
  \end{subfigure}
  \caption{Sample geometries for testing Green's identity,
    verified at targets on the boundary and in the exterior domain.}
\end{figure}

For the Helmholtz parameter $\omega=12.43$, let the boundary $\gamma$,
$\gamma(t) = (x_{1}(t),x_{2}(t))$, be parametrized as:
\[
   x_{1}(t) = \text{Re} \left( \sum_{j=0}^{50} \hat{x}_{1,j} \,
                            e^{2 \pi i j t} \right) ,\quad
   x_{2}(t) = \text{Re} \left( \sum_{j=0}^{50} \hat{x}_{2,j} \,
                            e^{2 \pi i j t} \right)  ,
\]
where the Fourier coefficients $\hat{x}_{1,j}, \hat{x}_{2,j}$ are
listed in~\ref{sec:coeffs}. This parametrization traces a fish-like
boundary, see Figure~\ref{fig:sample_geom0}.
Following the procedure described in \cite{qbx}, given a panel
order $q$, and a tolerance $\epsilon$, we refine $\gamma$ into
piecewise panels, with the functions $x_{1}$, $x_{2}$ on each panel
 interpolated using a $q$-term
Legendre  polynomial expansions. The panels are refined until the each
expansion and its requisite derivatives are resolved to the prescribed
tolerance $\epsilon$ in a spectral $\ell_2$-sense.

A typical test domain for all the numerical examples in this section is
described here.
Let $D_{1},\ldots,D_N$ be a collection of obstacles whose boundaries
$\Gamma_{j}$ are $\gamma$ up to an affine transformation, see
Figure~\ref{fig:sample_geom}.
Let $\Gamma = \cup_{j} \Gamma_{j}$, and
$\Omega^c= \mathbb{R}^{2} \setminus \cup_{j} D_{j}$ denote the exterior of these
obstacles.

Suppose $u$ satisfies the Helmholtz equation in $\Omega^{c}$ along with
the Sommerfeld radiation condition given by equations
\eqref{eq:helm-eq}, and \eqref{eq:sommerfeld-rad} respectively;
then $u$ satisfies the Green's identity
\begin{equation}
  u = \mathcal{D}\left[ u \right] -
                \mathcal{S} \left[ \frac{\partial u}{\partial n} \right]
  \label{eq:greensid}
\end{equation}
everywhere in $\Omega^c$. We verify this identity using QBX to
evaluate the layer potentials $\cS$ and $\cD$ arbitrarily close to
the boundary.

\subsection{Accuracy and Complexity \label{sec:acc-comp}}
Let $u$ be the Helmholtz potential generated by point sources
placed inside the
domains $D_{j}$. The potential is then given by
\[
u = \sum_{j=1}^{N} q_{j} \, H^{(1)}_0 (\omega \vert \bx- \bx_{j} \vert),
\]
where $\bx_{j} \in D_{j}$, and $q_{j}$ are randomly chosen.
Obviously $u$ satisfies the homogeneous Helmholtz equation
in $\Omega^c$,  and therefore satisfies identity~\eqref{eq:greensid}.
To test the accuracy of the algorithm in~\ref{sec:complete-fmm},
we compute the layer potentials
$\mathcal{S}[ \partial u/\partial n]$, $\mathcal{D} [ u ]$ and
obtain the error in \eqref{eq:greensid} at targets on the
boundary~$\Gamma$ and in the exterior $\Omega^c$.
Note that on the boundary, relationship~\eqref{eq:greensid} is
interpreted in the one-sided limit as in~\eqref{eq_oneside}.
The targets in $\Gamma_{\text{near}}$ are identified using the algorithm
described in Section~\ref{target-flagging}.

\begin{rem}
For all of the numerical experiments, we use a {\em level-restricted}
quad-tree for sorting sources, targets, and expansion centers in
the computational domain.
In a level-restricted quad-tree, two childless
boxes which share a boundary point are no more than one level of
refinement apart.
There are several standard algorithms for converting
a fully adaptive quad-tree into a level-restricted quad-tree
and we implement the
one discussed in~\cite{ethridge2000fast}.
\end{rem}

Let $\epsilon_{u,b}$ and $\epsilon_{u,v}$ be the weighted $\ell_2$
relative-error in Green's identity at the targets on the boundary and
the volume, respectively:
\[
\epsilon^2_{u,b} = \frac{\sum_{j=1}^{n_{t}} \left| u(\bt_{j}) -
u_{\text{qbx}} (\bt_{j}) \right|^{2} \, w_{j}}{ \sum_{j=1}^{n_{t}}
\left|u ( \bt_{j} ) \right|^{2} \, w_{j}}  ,
\]
and
\[
\epsilon^2_{u,v} = \frac{\sum_{j=1}^{n_{t}} \left| u(\bt_{j}) -
u_{\text{qbx}} (\bt_{j}) \right|^{2}}{ \sum_{j=1}^{n_{t}}
\left|u ( \bt_{j} ) \right|^{2}} .
\]
Here $u_{\text{qbx}} = \mathcal{S} [ {\partial u}/{\partial n}]
- \mathcal{D} [u]$ is computed using QBX and
$w_{j}$ is the Gaussian quadrature weight at the corresponding source
on the boundary.
Thus $\epsilon_{u,b}$ is a numerical approximation to the continuous
relative $L^{2}$ error on the boundary:
\[
\frac{\int_{\Gamma} |u-u_{\text{qbx}}
|^2 ds}{\int_{\Gamma} |u|^2 ds}.
\]

\begin{figure}[t]
\begin{center}
\includegraphics[width=5cm]{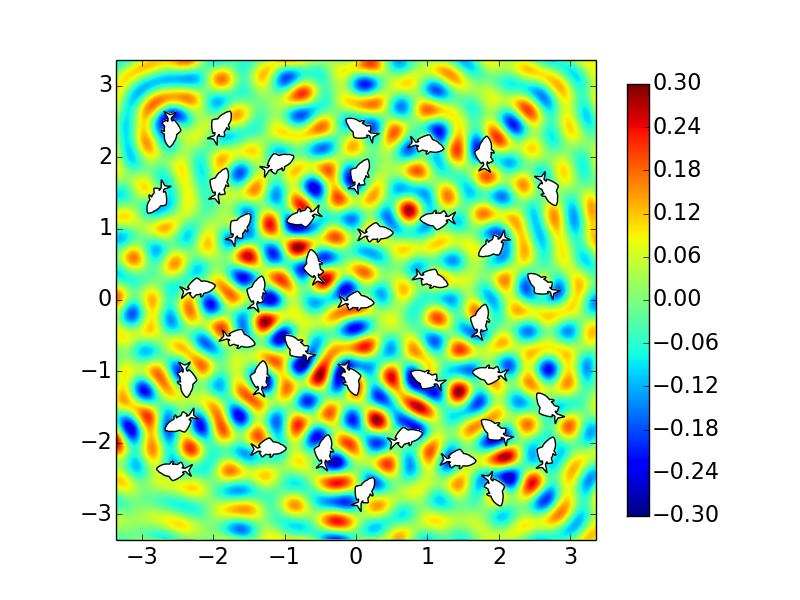}
\includegraphics[width=5cm]{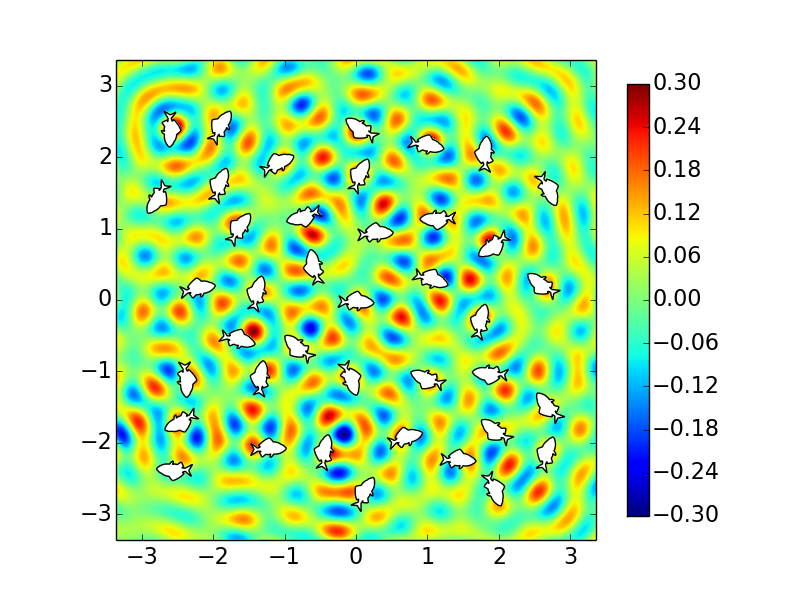}
\includegraphics[width=5cm]{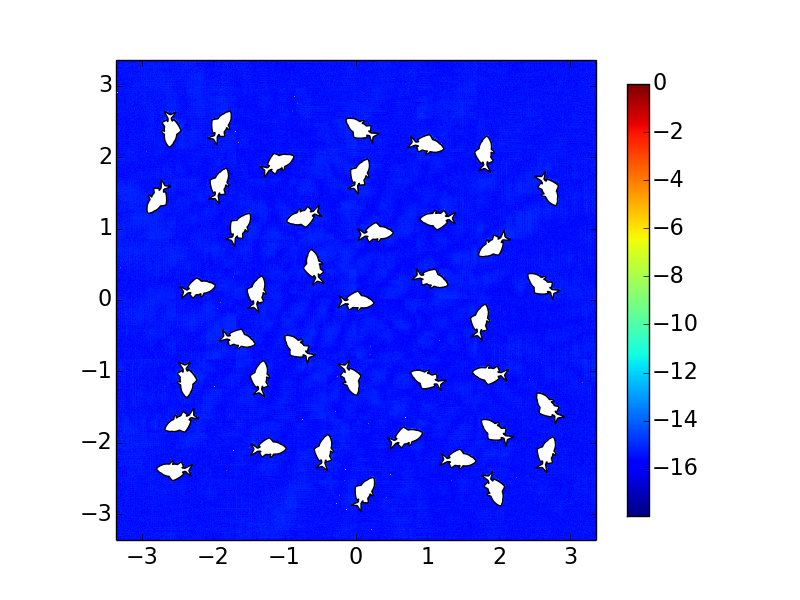}
\end{center}
\caption{Contour plots for (left) real part of $u$, (middle) imaginary part
of $u$, and (right) error in Green's identity. Boundary discretized using
order 16 Gauss-Legendre panels with 304640 points and potential evaluated at
600166 targets in the volume.}
\label{fig:greentest_samp}
\end{figure}

We also analyze the performance of the algorithm for different combinations
of $q$, $p$, and~$\epsilon$.
The results are summarized in Tables~\ref{tab-boundary}
and~\ref{tab-volume}.
The first column $\epsilon$ is the tolerance requested in the algorithm.
The second column is the order of Gauss-Legendre panels, given by $q$.
The third column is the QBX expansion order $p$.
Columns 4-7, denoted by $n_{d}$, $n_{s}$, $n_{t}$, $n_{e}$, are the number of
discretization nodes on $\Gamma$,
the number of over-sampled nodes on $\Gamma$, the number of targets,
and the number of expansion centers, respectively.
Columns 8-9, denoted
by $\epsilon_{s}$, $\epsilon_{d}$, are the resolution of the single-layer
 density $\sigma$, and the double-layer density $\mu$, respectively.
Let $a_{j,k}$ denote the coefficients of the
Legendre expansion of a function $f$ on panel $\Gamma_{k}$.
The resolution of the function $f$, denoted $\epsilon_{f}$,
on the discretization of the geometry is then given by
\begin{equation}
    \epsilon^2_{f} = \max_{k} \, \frac{
                        \sum_{j=q-n_{\text{tail}} }^{q} |a_{j,k}|^2}
                  {\sum_{j=1}^{q} |a_{j,k}|^2} \, h_{k} \, ,
\end{equation}
where $n_{\text{tail}}=1,2,$ or $3$ depending on the panel order $q$.
The error $\epsilon_{f}$ is the maximum relative $\ell^2$-norm of the tail
of the Legendre expansion of $f$ scaled by the arclength of
the panel.
Column $10$ is $\epsilon_{u,b}$ for Table~\ref{tab-boundary}
and~$\epsilon_{u,v}$ for Table~\ref{tab-volume}.
Finally, columns 11-13, denoted by $t_{qbx}$, $t_{fmm,1}$, $t_{fmm,2}$,
are the computation times.
The time $t_{qbx}$ is the time required to evaluate the layer potential,
$t_{fmm,1}$ is the time required for an $FMM$ with $n_{s}$
sources and $n_{t}$ targets, and $t_{fmm,2}$ is the time required
for an FMM with $n_{d}$ sources and $n_{t}$ targets. A plot of the
potential and errors in the Green's identity test is given in
Figure~\ref{fig:greentest_samp}.

\begin{rem}
The difference between $t_{qbx}$ and $t_{fmm,1}$ is exactly the additional
computational work in the FMM-accelerated QBX algorithm over the
original FMM.
The time $t_{fmm,2}$, on the other hand, is the time required to apply an
FMM on the original distribution of sources and targets.
\end{rem}

\begin{table}[H]
{\small
\centering
\resizebox{\columnwidth}{!}{%
\begin{tabular}{|c|c|c|c|c|c|c|c|c|c|c|c|c|}
\hline
$\epsilon$ & $q$ & $p$ & $n_{d}$ & $n_{s}$ & $n_{t}$ & $n_{e}$ & $\epsilon_{s}$ & $\epsilon_{d}$ & $\epsilon_{u}$ & $t_{qbx}$ & $t_{fmm,1}$ & $t_{fmm,2}$ \\ \hline
\multirow{4}{*} {5.00e-04} & \multirow{4}{*} {2} & \multirow{4}{*} {2} & 82536 & 330144 & 82536 & 82536 & 1.81e-04 & 7.49e-06 & 1.43e-03 & 11.71 & 4.24 & 2.00 \\ \hhline{~~~----------}
 & & & 165072 & 660288 & 165072 & 165072 & 6.85e-05 & 6.06e-06 & 1.53e-03 & 23.90 & 8.29 & 3.91 \\ \hhline{~~~----------}
 & & & 247608 & 990432 & 247608 & 247608 & 6.58e-05 & 1.75e-05 & 1.54e-03 & 37.65 & 12.78 & 6.02 \\ \hhline{~~~----------}
 & & & 330144 & 1320576 & 330144 & 330144 & 1.07e-04 & 1.04e-04 & 1.45e-03 & 48.67 & 17.04 & 8.07 \\ \hline
\multirow{4}{*} {5.00e-07} & \multirow{4}{*} {4} & \multirow{4}{*} {4} & 59328 & 296640 & 59328 & 59328 & 4.13e-06 & 1.52e-07 & 3.01e-07 & 15.28 & 5.73 & 2.27 \\ \hhline{~~~----------}
 & & & 118656 & 593280 & 118656 & 118656 & 2.28e-06 & 2.06e-07 & 3.06e-07 & 32.01 & 12.05 & 4.83 \\ \hhline{~~~----------}
 & & & 177984 & 889920 & 177984 & 177984 & 1.03e-06 & 6.69e-08 & 3.07e-07 & 46.87 & 17.48 & 7.26 \\ \hhline{~~~----------}
 & & & 237312 & 1186560 & 237312 & 237312 & 5.36e-06 & 1.90e-06 & 3.17e-07 & 61.44 & 22.00 & 9.10 \\ \hline
\multirow{4}{*} {5.00e-10} & \multirow{4}{*} {8} & \multirow{4}{*} {6} & 68544 & 342720 & 68544 & 68544 & 1.93e-10 & 9.49e-13 & 2.01e-10 & 36.38 & 9.35 & 3.92 \\ \hhline{~~~----------}
 & & & 137088 & 685440 & 137088 & 137088 & 3.24e-10 & 9.62e-13 & 6.22e-11 & 70.60 & 18.66 & 7.81 \\ \hhline{~~~----------}
 & & & 205632 & 1028160 & 205632 & 205632 & 2.91e-09 & 1.89e-12 & 1.78e-10 & 105.01 & 27.44 & 11.75 \\ \hhline{~~~----------}
 & & & 274176 & 1370880 & 274176 & 274176 & 3.51e-10 & 1.07e-12 & 1.00e-10 & 146.54 & 37.96 & 15.75 \\ \hline
\multirow{4}{*} {5.00e-13} & \multirow{4}{*} {16} & \multirow{4}{*} {8} & 76160 & 304640 & 76160 & 76160 & 4.19e-15 & 3.45e-16 & 3.79e-12 & 61.37 & 14.15 & 5.23 \\ \hhline{~~~----------}
 & & & 152320 & 609280 & 152320 & 152320 & 6.72e-15 & 1.03e-15 & 5.00e-12 & 118.22 & 26.96 & 10.80 \\ \hhline{~~~----------}
 & & & 228480 & 913920 & 228480 & 228480 & 1.13e-14 & 2.94e-15 & 4.02e-12 & 175.58 & 40.15 & 13.14 \\ \hhline{~~~----------}
 & & & 304640 & 1218560 & 304640 & 304640 & 8.76e-15 & 1.61e-15 & 6.36e-12 & 226.35 & 51.57 & 18.18 \\ \hline
\end{tabular}
}
\caption{Targets on the boundary}
\label{tab-boundary}
}
\end{table}

\begin{table}[H]
{\small
\centering
\resizebox{\columnwidth}{!}{%
\begin{tabular}{|c|c|c|c|c|c|c|c|c|c|c|c|c|}
\hline
$\epsilon$ & $q$ & $p$ & $n_{d}$ & $n_{s}$ & $n_{t}$ & $n_{e}$ & $\epsilon_{s}$ & $\epsilon_{d}$ & $\epsilon_{u}$ & $t_{qbx}$ & $t_{fmm,1}$ & $t_{fmm,2}$ \\ \hline
\multirow{4}{*} {5.00e-04} & \multirow{4}{*} {2} & \multirow{4}{*} {2} & 55024 & 220096 & 550240 & 55024 & 1.24e-04 & 1.49e-06 & 2.13e-05 & 9.48 & 3.38 & 2.28 \\ \hhline{~~~----------}
 & & & 55024 & 220096 & 1100480 & 55024 & 1.24e-04 & 1.49e-06 & 1.74e-05 & 16.41 & 5.09 & 3.97 \\ \hhline{~~~----------}
 & & & 55024 & 220096 & 1650720 & 55024 & 1.24e-04 & 1.49e-06 & 1.61e-05 & 17.42 & 13.72 & 12.44 \\ \hhline{~~~----------}
 & & & 55024 & 220096 & 2200960 & 55024 & 1.24e-04 & 1.49e-06 & 1.68e-05 & 16.33 & 8.73 & 7.53 \\ \hline
\multirow{4}{*} {5.00e-07} & \multirow{4}{*} {4} & \multirow{4}{*} {4} & 59328 & 296640 & 593280 & 59328 & 4.13e-06 & 1.52e-07 & 1.99e-09 & 18.55 & 9.06 & 6.32 \\ \hhline{~~~----------}
 & & & 59328 & 296640 & 1186560 & 59328 & 4.13e-06 & 1.52e-07 & 8.53e-09 & 21.00 & 8.20 & 5.15 \\ \hhline{~~~----------}
 & & & 59328 & 296640 & 1779840 & 59328 & 4.13e-06 & 1.52e-07 & 2.41e-09 & 35.51 & 14.55 & 11.48 \\ \hhline{~~~----------}
 & & & 59328 & 296640 & 2373120 & 59328 & 4.13e-06 & 1.52e-07 & 2.69e-09 & 33.04 & 25.13 & 22.74 \\ \hline
\multirow{4}{*} {5.00e-10} & \multirow{4}{*} {8} & \multirow{4}{*} {6} & 58752 & 293760 & 587520 & 58752 & 1.49e-10 & 5.15e-13 & 1.13e-13 & 40.46 & 9.70 & 5.12 \\ \hhline{~~~----------}
 & & & 58752 & 293760 & 1175040 & 58752 & 1.49e-10 & 5.15e-13 & 1.06e-13 & 41.40 & 12.76 & 8.09 \\ \hhline{~~~----------}
 & & & 58752 & 293760 & 1762560 & 58752 & 1.49e-10 & 5.15e-13 & 9.56e-14 & 46.17 & 14.60 & 9.80 \\ \hhline{~~~----------}
 & & & 58752 & 293760 & 2350080 & 58752 & 1.49e-10 & 5.15e-13 & 1.19e-13 & 74.95 & 18.56 & 13.60 \\ \hline
\multirow{4}{*} {5.00e-13} & \multirow{4}{*} {16} & \multirow{4}{*} {8} & 60928 & 243712 & 609280 & 60928 & 4.80e-15 & 1.91e-16 & 1.20e-14 & 52.95 & 15.60 & 8.86 \\ \hhline{~~~----------}
 & & & 60928 & 243712 & 1218560 & 60928 & 4.80e-15 & 1.91e-16 & 3.21e-14 & 65.63 & 20.04 & 13.22 \\ \hhline{~~~----------}
 & & & 60928 & 243712 & 1827840 & 60928 & 4.80e-15 & 1.91e-16 & 2.74e-14 & 71.13 & 23.02 & 15.09 \\ \hhline{~~~----------}
 & & & 60928 & 243712 & 2437120 & 60928 & 4.80e-15 & 1.91e-16 & 2.45e-14 & 83.82 & 30.18 & 22.44 \\ \hline
\end{tabular}
}
\caption{Targets in the volume}
\label{tab-volume}
}
\end{table}



\subsection{A scattering problem}
{\em Sound-soft} scattering problems in acoustics correspond to exterior
Dirichlet boundary value problems for the Helmholtz equation.
Let $u_{tot}$, $u_{inc}$, $u_{sc}$ denote the total potential,
the incident potential and the
scattered potential, respectively, all of which solve the homogeneous
Helmholtz equation in
the exterior region of a collection of obstacles
except possibly at a finite number of points.
Given an incident potential $u_{inc}$, the goal is to compute the
scattered potential $u_{sc}$ such that $u_{tot} = 0$ on the boundary
$\Gamma$, where $\Gamma$ denotes the boundary of the obstacles.
Thus $u_{sc}$ solves the exterior Helmholtz Dirichlet problem
given by equations \eqref{eq:helm-eq}, \eqref{eq_diri}, and
\eqref{eq:sommerfeld-rad}, with Dirichlet data
$f = -u_{inc}$ on the boundary $\Gamma$.

\begin{figure}[t]
\begin{center}
\includegraphics[width=6cm]{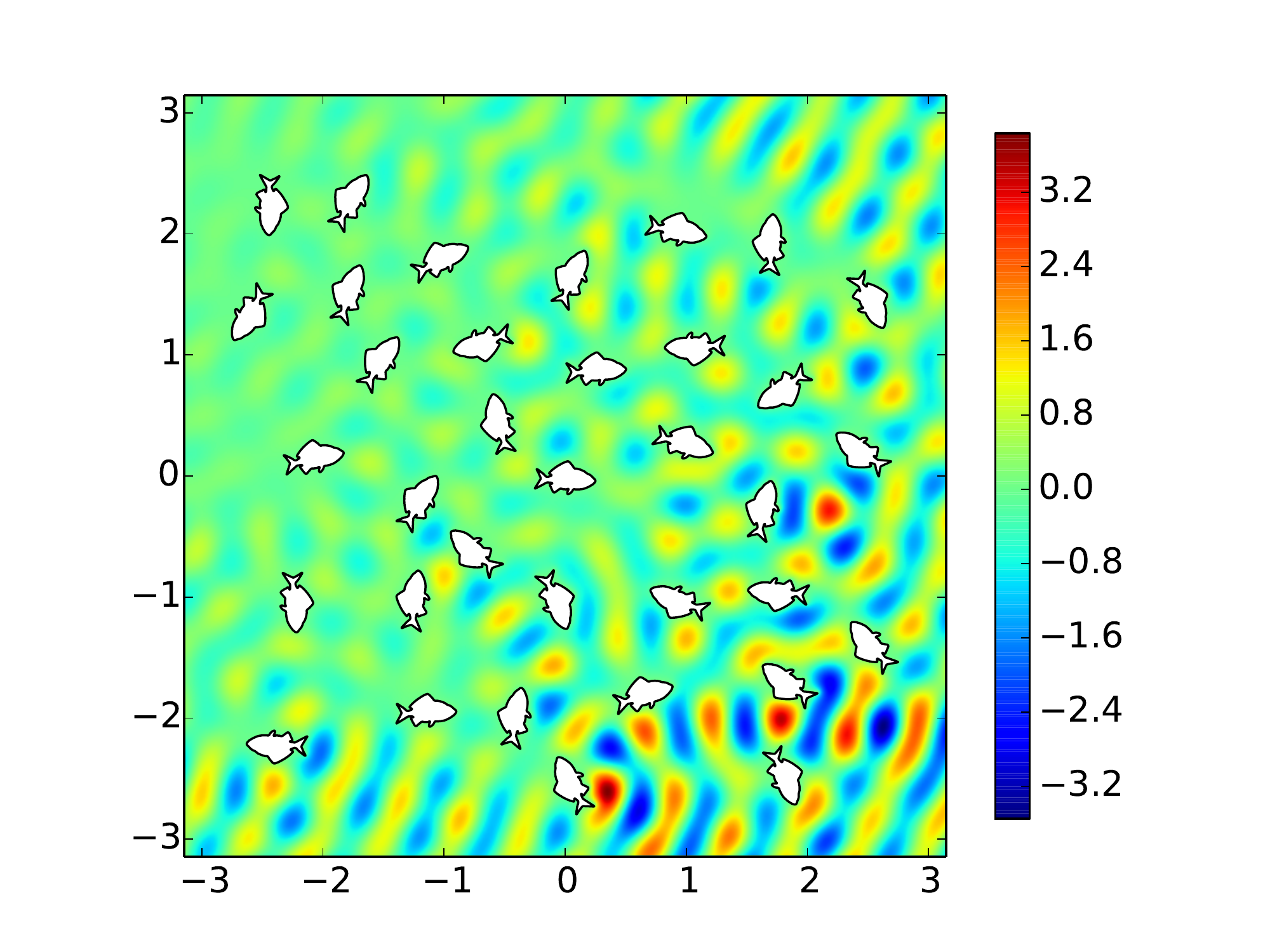}
\includegraphics[width=6cm]{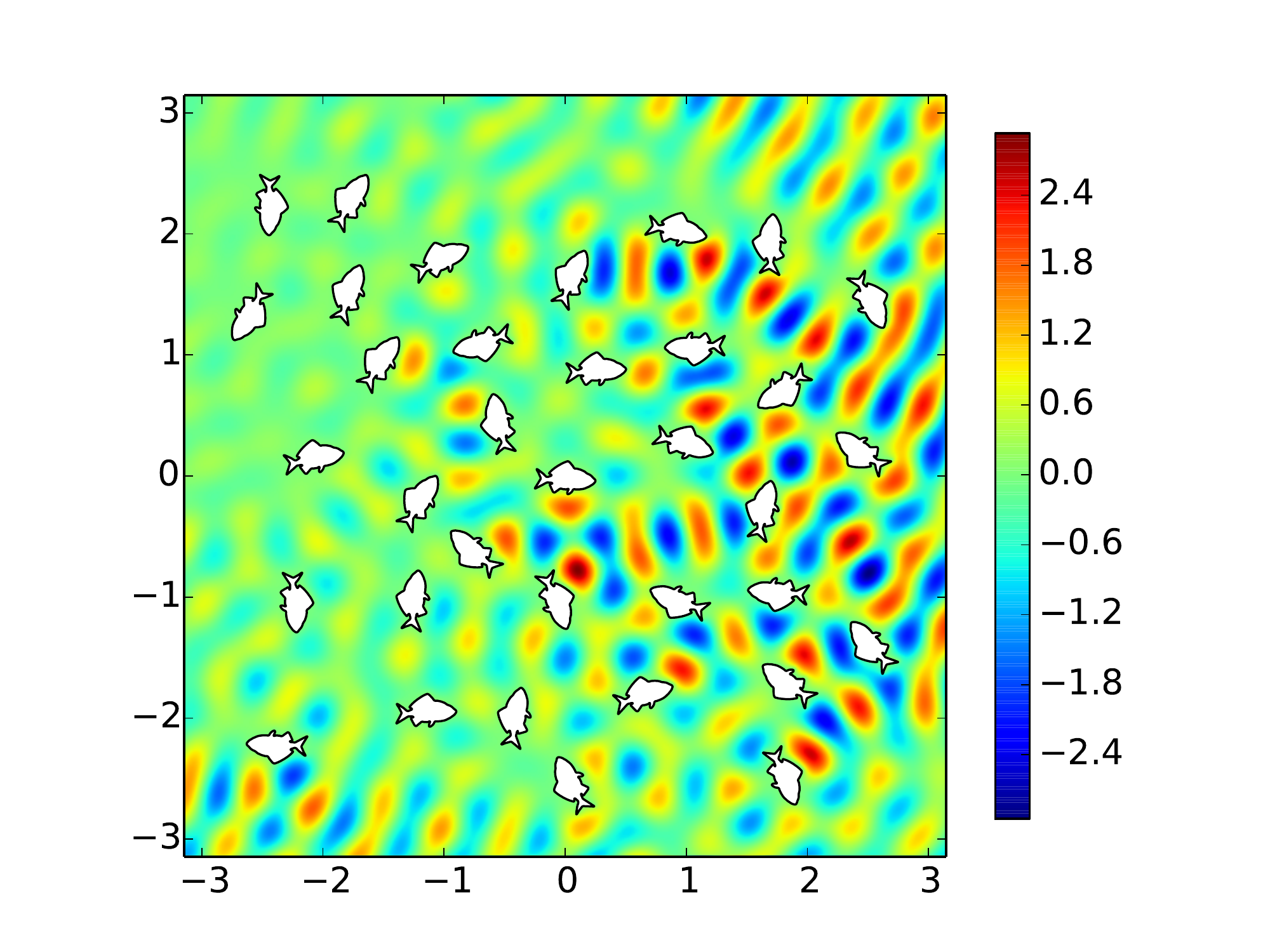}
\end{center}
\caption{Contour plots for (left) real part of $u_{tot}$, (right) imaginary
part of $u_{tot}$ evaluated at $600 113$ target points in the
volume. The target points are allowed to be arbitrarily close to the objects.}
\label{fig-cont}
\end{figure}

Following the procedure described in the introduction, we represent the
scattered potential as
\begin{equation}
\label{eq_rep}
u_{sc} = \mathcal{D}[ \sigma  ] + i\omega \mathcal{S}
[ \sigma ] ,
\end{equation}
where $\sigma$ is an unknown density.
On imposing the boundary conditions and using properties of the single- and
double-layer potentials, we have the integral equation along $\Gamma$:
\begin{equation}
\label{eq-int_eq}
\frac{1}{2}\sigma + \mathcal{D}^{*} [\sigma ] + i\omega
        \mathcal{S}^{*} [\sigma] = -u_{inc}.
\end{equation}
We discretize the above equation using a Nystr\"{o}m method and use
Global QBX for computing the layer potentials.
Using the geometry discretization described in
Section~\ref{sec-num-prel}, consider
the scattering problem in the exterior of $35$ inclusions discretized with
$16^\text{th}$-order Legendre expansions ($q =16$)
 and $\epsilon = 5.0\times 10^{-7}$.
Each panel is further subdivided once to ensure that the solution $\sigma$ is
well-resolved.
Let $u_{inc}$ be a plane-wave given by
\begin{equation}
  u_{inc} (\bx) = e^{i \omega \left(-2x_{1} + x_{2} \right)/\sqrt{5}} ,
\end{equation}
where $\bx = \left(x_{1},x_{2} \right)$.
On discretizing integral equation~\eqref{eq-int_eq}, we have a linear
system with $105 280$ unknowns.
We use an iterative GMRES-based solver to obtain the solution $\sigma$;
iterations are performed until we reach a relative residual
of $1.0 \times 10^{-5}$.
The solution converged in $554$ iterations,
the resolution of the density
$\sigma$ on the given discretization was $9.41\times 10^{-6}$,
and the time required per iteration was \mbox{$31$ sec}.
In Figure~\ref{fig-cont}, we plot the real and imaginary parts of the
resulting total potential $u_{tot}$.

\section{Conclusions and future work}
\label{sec:conclusions}

We have introduced a method for scalably performing singular
quadrature using Quadrature by Expansion (QBX) within a fast algorithm
based on the fast multipole method (FMM). The resulting algorithm is
known as {\em FMM-Accelerated Global QBX}.

We have demonstrated that that the globally valid expansions of layer
potentials necessary to carry out QBX can scalably be constructed
using standard translation operators applied to $H$- and
$J$-expansions.
The resulting scheme scales like the underlying FMM in the low-frequency
regime, with an asymptotic  run time of $\mathcal O(n)$, where $n$
is proportional to the number of
sources in the discretization and targets in the volume. The constant
implicit in the $\mathcal O(\cdot)$ notation is only a small factor
larger than that in the point-based FMM, roughly between two and four,
depending on the desired precision.

Beyond merely providing a method for the evaluation of layer
potentials, the scheme not only verifies a number of conditions
required to guarantee its accuracy, it also automatically aids the user
(by ways of mesh refinement) in ensuring that these conditions are
met.
We have thus obtained a capability to accurately evaluate layer
potentials, scalably, anywhere in space, in a black-box fashion.

Under some circumstances, it is possible to construct algorithms that
require less mesh refinement than the method presented here. This is
particularly true in the case when some parts of the source
geometry come so close to other parts that they \emph{almost} touch
each other. While
the presented algorithm \emph{will} provide an accurate answer case,
it may not do so with the best possible efficiency.
Motivated by these perspectives and the promising performance
of the method presented here, a scheme denoted
as {\em Local QBX} will be introduced in a subsequent
paper~\cite{localqbx} currently in preparation.
This alternative scheme constructs expansions of the potential due
to smaller pieces of the geometry, as necessary, in order to overcome
extra geometry refinement.

Lastly, it is relatively straightforward to derive Global QBX schemes
for computing layer potentials due to Stokeslets and stresslets in
fluid dynamics, current and charge densities in electromagnetics, as
well as other
classical potentials in mathematical physics~\cite{rahimian_2016}.
These, and extensions to
three dimensions, are ongoing projects.


\section*{Acknowledgements}
M.~Rachh's research was partially supported by the U.S.
Department of Energy under contract DEFG0288ER25053 and by the Office
of the Assistant Secretary of Defense for Research and Engineering and
AFOSR under NSSEFF Program Award FA9550-10-1-0180.
A. Klöckner's research was supported by the
National Science Foundation under grant DMS-1418961 and by the Office
of the Assistant Secretary of Defense for Research and Engineering and
AFOSR under NSSEFF Program Award FA9550-10-1-0180.
M. O'Neil's research was
supported in part by the Office of the Assistant Secretary of
Defense for Research and Engineering and AFOSR under NSSEFF
Program Award FA9550-10-1-0180, and by the Office of Naval Research
under Award N00014-15-1-2669.
The authors would like to thank Leslie Greengard, Abtin Rahimian, and
Matt Wala for several useful discussions.

The authors would further thank the anonymous referees for many helpful comments
that led to a much-improved manuscript.

\appendix
\section{Fourier coefficients for test geometry}
\label{sec:coeffs}
\begin{table}[H]
{\small
\centering
\begin{tabular}{|c|c|c|}
\hline
$j$ & $\hat{x}_{1,j}$ & $\hat{x}_{2,j}$ \\ \hline
0 & -3.03e-02 + i 0.00e+00 & -1.56e-02 + i 0.00e+00 \\ \hline
1 & -1.00e-01 + i 2.34e-02 & -1.01e-02 - i 4.92e-02 \\ \hline
2 & 1.28e-02 + i 2.16e-03 & -1.50e-03 + i 1.37e-02 \\ \hline
3 & -9.40e-03 + i 3.98e-03 & 2.23e-04 - i 5.08e-03 \\ \hline
4 & 3.18e-03 - i 1.92e-03 & -7.06e-03 - i 2.70e-03 \\ \hline
5 & -3.42e-03 + i 3.37e-03 & -9.79e-03 + i 2.51e-03 \\ \hline
6 & -2.13e-03 - i 3.87e-03 & -3.70e-03 + i 4.34e-03 \\ \hline
7 & -4.24e-03 - i 3.45e-03 & -2.36e-03 - i 1.83e-03 \\ \hline
8 & -1.61e-03 - i 2.24e-03 & 2.46e-03 - i 8.88e-04 \\ \hline
9 & -1.32e-03 + i 1.85e-03 & 2.16e-03 - i 1.10e-04 \\ \hline
10 & -4.06e-04 + i 2.52e-04 & -1.92e-03 - i 6.47e-04 \\ \hline
11 & 5.58e-04 + i 6.41e-04 & -1.07e-03 - i 1.74e-04 \\ \hline
12 & -1.29e-04 + i 1.88e-04 & -4.82e-05 - i 1.37e-04 \\ \hline
13 & -8.71e-04 + i 1.47e-03 & -1.60e-03 - i 4.67e-05 \\ \hline
14 & -5.12e-04 - i 1.51e-04 & 4.71e-04 - i 3.93e-04 \\ \hline
15 & 4.31e-04 - i 4.82e-04 & 6.10e-04 + i 2.67e-04 \\ \hline
16 & -3.51e-04 - i 5.28e-04 & -5.65e-04 + i 8.88e-04 \\ \hline
17 & -7.72e-04 - i 2.93e-04 & -5.04e-04 - i 2.83e-04 \\ \hline
18 & -3.65e-04 - i 2.33e-04 & 1.37e-04 - i 4.91e-04 \\ \hline
19 & 8.68e-04 + i 3.97e-04 & 9.03e-05 + i 9.22e-05 \\ \hline
20 & 1.50e-04 + i 1.72e-04 & -2.04e-04 - i 4.82e-05 \\ \hline
21 & -2.22e-04 - i 1.72e-04 & -3.64e-04 - i 2.16e-04 \\ \hline
22 & -3.09e-04 + i 2.04e-05 & -3.61e-04 + i 3.53e-05 \\ \hline
23 & -1.92e-04 + i 2.53e-04 & -1.04e-04 - i 1.73e-05 \\ \hline
24 & -3.36e-04 - i 1.48e-04 & 7.34e-05 + i 1.40e-04 \\ \hline
25 & -1.16e-04 - i 6.38e-04 & 8.94e-05 - i 1.08e-04 \\ \hline
\end{tabular}
\hspace{1ex}
\begin{tabular}{|c|c|c|}
\hline
$j$ & $\hat{x}_{1,j}$ & $\hat{x}_{2,j}$ \\ \hline
26 & -2.24e-04 - i 1.73e-04 & 8.66e-05 - i 2.07e-05 \\ \hline
27 & -5.47e-05 + i 2.16e-04 & 3.60e-05 - i 1.05e-04 \\ \hline
28 & 9.47e-05 + i 3.04e-04 & -3.38e-04 + i 4.25e-06 \\ \hline
29 & 1.87e-04 + i 7.48e-05 & -5.85e-05 - i 7.12e-05 \\ \hline
30 & -6.42e-05 + i 2.08e-05 & -1.01e-04 - i 4.42e-05 \\ \hline
31 & -2.33e-04 + i 2.49e-05 & -3.08e-05 + i 6.74e-05 \\ \hline
32 & -1.47e-04 + i 7.06e-05 & 7.47e-05 + i 3.94e-05 \\ \hline
33 & 3.51e-05 - i 1.69e-04 & -3.73e-05 - i 4.19e-06 \\ \hline
34 & 4.50e-05 - i 1.88e-04 & -1.20e-04 + i 3.74e-05 \\ \hline
35 & -9.51e-05 - i 1.18e-04 & -1.00e-05 - i 7.77e-05 \\ \hline
36 & -8.54e-05 + i 7.05e-05 & -7.14e-05 - i 6.35e-05 \\ \hline
37 & 9.22e-05 + i 9.62e-05 & -1.57e-05 - i 7.09e-05 \\ \hline
38 & 1.07e-04 + i 5.55e-05 & 2.40e-05 - i 1.28e-04 \\ \hline
39 & -5.84e-05 - i 5.48e-05 & -8.74e-05 + i 1.17e-04 \\ \hline
40 & -1.58e-04 - i 4.45e-05 & -9.08e-05 + i 1.12e-05 \\ \hline
41 & -1.31e-04 - i 2.73e-05 & 3.18e-05 - i 4.73e-05 \\ \hline
42 & -6.19e-06 - i 2.10e-05 & 1.22e-04 + i 4.18e-05 \\ \hline
43 & -8.43e-06 - i 7.75e-05 & -2.89e-05 + i 3.22e-05 \\ \hline
44 & -5.35e-05 - i 2.64e-05 & -1.11e-04 - i 3.66e-05 \\ \hline
45 & -2.68e-06 + i 1.33e-05 & -3.82e-05 - i 6.75e-05 \\ \hline
46 & 4.99e-05 + i 1.14e-04 & -4.55e-05 - i 1.41e-05 \\ \hline
47 & 6.65e-06 + i 4.98e-05 & -2.83e-05 - i 5.38e-05 \\ \hline
48 & -2.05e-05 - i 6.93e-05 & -2.80e-05 - i 1.56e-05 \\ \hline
49 & -2.32e-05 - i 6.10e-05 & 2.21e-05 + i 1.22e-05 \\ \hline
50 & -2.31e-05 + i 3.32e-05 & 5.11e-05 + i 4.80e-05 \\ \hline
\end{tabular}
\caption{Fourier coefficients}
\label{tab-coeffs}
}
\end{table}

\section*{References}

\bibliography{qbx-accel}
\bibliographystyle{abbrvnat}

\end{document}